\documentclass[journal]{IEEEtran}
\usepackage{cite}

\ifCLASSINFOpdf
  \usepackage[pdftex]{graphicx}
  \graphicspath{{../Graphics/}}
\else
\fi
\usepackage[cmex10]{amsmath}
\usepackage{labtex}
\usepackage[caption=false,font=footnotesize]{subfig}

\usepackage{stfloats}
\usepackage{url}

\usepackage{graphicx}
\usepackage{fancyhdr}
\usepackage{verbatim}
\usepackage{listings}
\usepackage{amsmath}
\usepackage{amssymb}
\usepackage{sidecap}
\usepackage{array}
\usepackage[usenames, dvipsnames]{color}
\usepackage{pgfplots}
\pgfplotsset{compat=1.6}

\DeclareFontFamily{U}{msb}{}
\DeclareFontShape{U}{msb}{m}{n}{ <5> <6> <7> <8> <9> gen * msbm
<10> <10.95> <12> <14.4> <17.28> <20.74> <24.88> msbm10}{}
\DeclareSymbolFont{AMSb}{U}{msb}{m}{n}
\DeclareMathSymbol{\Reals}{\mathalpha}{AMSb}{'122}
\DeclareMathSymbol{\Naturals}{\mathalpha}{AMSb}{'116}
\DeclareMathSymbol{\Knumbers}{\mathalpha}{AMSb}{'113}
\DeclareMathSymbol{\Rationals}{\mathalpha}{AMSb}{'121}
\DeclareSymbolFont{AMSb}{U}{msb}{m}{n}
\DeclareMathSymbol{\setB}{\mathalpha}{AMSb}{'102}
\DeclareMathSymbol{\setC}{\mathalpha}{AMSb}{'103}
\DeclareMathSymbol{\setD}{\mathalpha}{AMSb}{'104}
\DeclareMathSymbol{\setE}{\mathalpha}{AMSb}{'105}
\DeclareMathSymbol{\setF}{\mathalpha}{AMSb}{'106}
\DeclareMathSymbol{\setI}{\mathalpha}{AMSb}{'111}
\DeclareMathSymbol{\setK}{\mathalpha}{AMSb}{'113}
\DeclareMathSymbol{\setM}{\mathalpha}{AMSb}{'115}
\DeclareMathSymbol{\setN}{\mathalpha}{AMSb}{'116}
\DeclareMathSymbol{\setP}{\mathalpha}{AMSb}{'120}
\DeclareMathSymbol{\setQ}{\mathalpha}{AMSb}{'121}
\DeclareMathSymbol{\setR}{\mathalpha}{AMSb}{'122}
\DeclareMathSymbol{\setS}{\mathalpha}{AMSb}{'123}

\usepackage{ifthen,version}
\newboolean{include-notes}
\setboolean{include-notes}{false}
\newcommand{\tmnote}[1]{\ifthenelse{\boolean{include-notes}}%
  {\textbf{(TM says: #1)}}{}}
\newcommand{\ejnote}[1]{\ifthenelse{\boolean{include-notes}}%
  {\textbf{(EJ says: #1)}}{}}
\newcommand{\jsnote}[1]{\ifthenelse{\boolean{include-notes}}%
  {\textbf{(JS says: #1)}}{}}

\newcommand{\deriv}[2][]{\dfrac{\partial #1}{\partial #2}}
\newcommand{\tderiv}[2][]{\tfrac{\partial #1}{\partial #2}}

\newcommand{\derivII}[3][]{\dfrac{\partial^2 #1}{\partial #2 \partial #3}}
\newcommand{\tderivII}[3][]{\tfrac{\partial^2 #1}{\partial #2 \partial #3}}

\newcommand{\Half}{\frac{1}{2}}
\newcommand{\half}{\tfrac{1}{2}}
\newcommand{\ud}{\mathrm{d}}

\newcommand{\op}{\circ}

\newcommand{\dq}{\ensuremath{\dot{q}}}

\newcommand{\Fp}{\ensuremath{F^+}}
\newcommand{\Fm}{\ensuremath{F^-}}
\newcommand{\K}{\mathcal{K}}

\newcommand{\dt}{\Delta t}

\providecommand{\e}[1]{\ensuremath{\times 10^{#1}}}

\newboolean{abbreviated-forms}
\setboolean{abbreviated-forms}{true}

\ifthenelse{\boolean{abbreviated-forms}}{
  \newcommand{\Ld}[1]{L_{#1}}
  \newcommand{\fdp}[1]{f_{#1}^+}
  \newcommand{\fdm}[1]{f_{#1}^-}
}{
  \newcommand{\Ld}[1]{
    \ifthenelse{#1=1}
               {L_d(q_{k-1},q_k)}
               {\ifthenelse{#1=2}
                 {L_d(q_k,q_{k+1})}
                 {\badargument}}}
  \newcommand{\fdp}[1]{
    \ifthenelse{#1=1}
               {f_d^+(q_{k-1},q_k,u_{k-1})}
               {\ifthenelse{#1=2}
                 {f_d^+(q_k,q_{k+1}, u_k)}
                 {\badargument}}}
  \newcommand{\fdm}[1]{
    \ifthenelse{#1=1}
               {f_d^-(q_{k-1},q_k, u_{k-1})}
               {\ifthenelse{#1=2}
                 {f_d^-(q_k,q_{k+1}, u_k)}
                 {\badargument}}}
}

\begin{document}
%
\title{Structured Linearization of Discrete Mechanical Systems for Analysis and
  Optimal Control}

\author{Elliot~Johnson, Jarvis~Schultz, and Todd~Murphey%
\thanks{Elliot~Johnson is with the Southwest Research Institute, San Antonio, TX, 78228 USA (e-mail:  elliot.johnson@swri.org).  His work was performed while a graduate student at Northwestern University.}
  \thanks{Jarvis~Schultz, and Todd~Murphey are with the Department of Mechanical Engineering,
    Northwestern University, Evanston, IL, 60201 USA. (e-mail:
     jschultz@u.northwestern.edu,
    t-murphey@northwestern.edu)}
}

%


\maketitle

\begin{abstract}
  Variational integrators are well-suited for simulation of mechanical systems
  because they preserve mechanical quantities about a system such as momentum,
  or its change if external forcing is involved, and holonomic
  constraints. While they are not energy-preserving they do exhibit long-time
  stable energy behavior. However, variational integrators often simulate
  mechanical system dynamics by solving an implicit difference equation at each
  time step, one that is moreover expressed purely in terms of configurations at
  different time steps.  This paper formulates the first- and second-order
  linearizations of a variational integrator in a manner that is amenable to
  control analysis and synthesis, creating a bridge between existing analysis
  and optimal control tools for discrete dynamic systems and variational
  integrators for mechanical systems in generalized coordinates with forcing and
  holonomic constraints.  The forced pendulum is used to illustrate the technique.  A second example solves the discrete LQR
  problem to find a locally stabilizing controller for a 40 DOF system with 6
  constraints.

  \noindent \textit{Note to Practitioners}---The practical value of
  this work is the explicit derivation of recursive formulas for exact
  expressions for the first- and second-order linearizations of an
  arbitrary constrained mechanical system without requiring symbolic
  calculations.  This is most applicable to the design of
  computer-aided design (CAD) software, where providing linearization
  information and sensitivity analysis facilitates mechanism analysis
  (e.g., controllability, observability) as well as control design
  (e.g., design of locally stabilizing feedback laws).
\end{abstract}

\begin{IEEEkeywords}
simulation, mechanism analysis, optimal control
\end{IEEEkeywords}

%
\IEEEpeerreviewmaketitle


\section{Introduction}

\IEEEPARstart{N}{umerical} integration schemes for mechanical systems typically
begin with continuous-time representations of dynamics (i.e., ordinary
differential equations) and then apply numerical integration to yield a
discrete-time approximation to the continuous-time dynamics. Instead of
computing the Euler-Lagrange equations based on extremizing the action integral,
variational integrators---sometimes referred to as \emph{structured integrators}
\cite{GeometricNumericalIntegration}---use a time-discretized form of the action
integral and then the resulting action sum is extremized.  The subsequent
integrators have the advantage of conserving momentum and a symplectic form as
well as bounding energy behavior \cite{lew2004b}.  Hence, variational
integrators avoid the calculation of any ordinary differential equations (ODE)
and lead to an implicit difference equation that is solved numerically at each
time step to find the next configuration of a mechanical
system.  
Additionally the numeric properties of variational integrators typically provide
stable energy behavior over long time horizons even with large timesteps
\cite{lew2004a} while also providing exact holonomic constraint
satisfaction\cite{LewThesis}.

Despite the fact that variational integrators simulate a system by numerically
solving an implicit non-linear equation at each time step, the process can be
abstractly represented in state form as an explicit, one-step discrete dynamic
system (i.e., $x_{k+1} = f(x_k,u_k)$) \cite{GeometricNumericalIntegration}.  The
contribution of this work is to calculate the first- and second-order
linearizations of this abstracted form, and it is useful that the linearization
may be explicitly calculated without ever calculating an analytic expression for the one-step map itself.  This
approach makes variational integrators compatible with a wide range of existing
analysis and optimal control tools for discrete dynamic systems (e.g., LQ
regulators \cite{AndersonMoore}, predictive control in assembly
\cite{Zhong2010}, deconvolution techniques \cite{You2011}).  The methods
described apply to complex mechanical systems in generalized coordinates with
many degrees of freedom and holonomic constraints, such as those previously
studied in \cite{JohnsonMurpheyTRO2008}.

The importance of linear systems analysis is evident in nearly all engineering
analysis, but computing linearizations can be analytically challenging (e.g.,
computing transverse linearizations for periodic systems \cite{Shiriaev10}) and
moreover computing the discrete time linearization from the continuous one is
susceptible to substantial numerical difficulties \cite{laub2010}.  By taking
advantage of discrete-time variational integrators we obtain exact
representations of the discrete-time linearization for arbitrary mechanical
systems subject to holonomic constraints.

Existing literature discusses linearizing variational integrators
\cite{HiraniThesis} for the purpose of constructing more accurate integrators.
That approach linearizes the discrete Lagrangian itself before the action
principle is applied.  The method described in this paper is specifically for
analysis and optimal control applications that require linearizations of the
dynamics, after the action principle is invoked.

Without the state form representation and linearizations discussed in this
paper, the implicit form of a variational integrator is incompatible with the
majority of existing analysis and optimal control techniques for discrete
systems described by an explicit first-order state form \cite{AndersonMoore}.
Applications that require more than simulation (e.g, stability analysis, design
of feedback control laws, optimal control) force the designer to either build an
auxiliary model with traditional continuous dynamics---that are only guaranteed
to apply to the discrete-time model in the limit as $dt\to0$---or develop new
methods specific to variational integrators.

The Discrete Mechanics Optimal Control (DMOC) framework
\cite{Ober2011,Leyendecker2010} offers such an approach to optimal control based
on variational integrators.  DMOC was developed to replace infinite dimensional,
continuous time optimal control problems with an (approximately) equivalent
finite-dimensional, constrained optimal control problem.  The optimal control
problem can then be solved using standard numeric optimization methods without
the problems associated with numeric integration and infinite dimensional vector
spaces.  One of the contributions of this paper is that we supplement DMOC by
connecting variational integrators with existing discrete optimal control
methods, providing the capability to take an optimal trajectory generated by
DMOC and generate a feedback regulator for it.  In \cite{sergio_discrete_2011}
DMOC is used to generate a discrete reference trajectory for the swing-up of a
cart-pendulum system.  In order to generate stabilizing feedback controllers
about this trajectory, the authors utilize interpolation to provide a continuous
representation of the discrete trajectory and thus allow standard
linearizations. A gain-scheduling technique is then used to piece together
solutions to a set of continuous LQR controllers, which are then presumably
implemented experimentally in discrete time.  The discrete linearizations
presented herein allow this entire process to occur in discrete time.
Similarly, the techniques discussed here could be used to generate feedback
laws for planning algorithms that provide plans in configuration space
\cite{ompl2012}.

In \cite{Ober2011}, it is shown that the discrete-time adjoint equation---the
equation that governs optimality of a control signal---for an explicit/implicit
partitioned Runge-Kutta scheme is itself an explicit/implicit partitioned
Runge-Kutta scheme.  Given that result, it is perhaps surprising that the
linearization of the (typically) implicitly defined variational integrator is in
fact an explicit calculation if the ``state'' of the system is chosen
appropriately.
It is by no means obvious that an implicit equation expressed directly
in terms of the configuration even has a linearization to which
classical methods in discrete-time optimal control can be applied.
The key observation is that because a variational integrator can be
rewritten as a one-step method (due to the existence of a so-called
generating function that generates the method and guarantees its
symplecticity), the one-step method provides the appropriate object to
linearize.  Moreover, another consequence of the existence of the
generating function is the existence of a modified Hamiltonian for the
system---a system of which the variational integrator is exactly
sampling at time steps $dt$.  This provides another interpretation of
the linearization we calculate here---it is both the \emph{exact}
discrete-time linearization of the variational integrator as well as
the \emph{exact} continuous-time state transition matrix for the
modified system evaluated at times that are $dt$ apart.  This
correspondence between the discrete-time interpretation and the
continuous-time interpretation indicates that the linearization of the
map is structure-preserving as well (e.g., the linearization does not
artificially introduce non-mechanical behavior by virtue of sampling
the trajectory).

The contribution of this paper is thus two-fold.  First, we show that
the linearization of a variational integrator---in this case that
obtained by using the midpoint rule---is explicit despite the
variational integrator being implicitly defined.  Moreover,
calculating the linearization is simply a matter of following a tree
structure that describes the mechanical system, even for constrained
systems.  We show the second derivative of a trajectory may be
obtained as well.  The consequence of this result is that
linearization information may be obtained in a purley algorithmic
foundation, without any symbolic computation; hence, it is reasonable
to expect computer aided design (CAD) software packages to include
linearization capability for arbitrary topologies for any given
numerical method to facilitate analysis and control design for complex
mechanisms, potentially operating in  scenarios where managing
sensitivity is crucial (see, for example, the editorial on this topic
\cite{dandrea2012} where software support of design is cited as a
major need in automation).  Secondly, we demonstrate that the control
calculation is well-posed even for high degree-of-freedom
systems, using a stabilization problem for a mechanical model of a 40
DOF marionette as an example.  We additionally briefly describe the
software package, \texttt{trep}, that we have written that implements
these techniques as well as corresponding techniques for
continuous-time dynamics (the continuous-time dynamics are discussed
in \cite{JohnsonACC2010}).

This paper is organized as follows.
Section~\ref{sec:variational-integrators} describes the particular
variational integrator used through the paper, and introduces a
pendulum example that is used to demonstrate the methods as they are
derived.  Section~\ref{sec:state-form} introduces the abstract
representation of a variational integrator as a first-order discrete
dynamic system.  After the background in
Secs.~\ref{sec:variational-integrators} and \ref{sec:state-form}, the
first- and second-order linearizations are derived in
Sec.~\ref{sec:first-deriv} and \ref{sec:second-deriv}, respectively.
The linearizations are extended to include systems with holonomic
constraints in Sec.~\ref{sec:constrained-systems}.  Section
\ref{sec:sing-linearization} presents several examples illustrating
singularities of the linearizations. An open source software
implementation called {\tt trep} is introduced in
Sec.~\ref{sec:implementation} and used to find stabilizing feedback
controllers first for the simple pendulum in
Sec.~\ref{sec:example:pend:optimization} and then for a humanoid
marionette in Sec.~\ref{sec:marionette} (this example serves as our
canonical example of a ``complex'' underactuated mechanism
\cite{egerstedthscc2007}).  Finally, Sec.~\ref{sec:conclusion}
summarizes the method, discusses the advantages and limitations, and
discusses future work.


\section{Variational Integrators \label{sec:variational-integrators}}

In this section we provide a brief overview of variational integrators
and present the specific variational integrator used throughout this paper.
More detailed introductions and discussions can be found in
\cite{Kharevych2006,WestThesis,LewThesis,MarsdenWest2001}.

The idea behind variational integrators is to discretize the action
with respect to time before finding the discrete-time equations of
motion.  Doing so leads to integration schemes that avoid problems
associated with numerically integrating a continuous ODE.  These
problems can occur because the numerical approximations that are
introduced do not respect fundamental mechanical properties like
conservation of momentum, energy, and a symplectic form, all of which
are relevant to mechanical systems (both forced and unforced).

The continuous-time dynamics of a mechanical system are described by
the Euler-Lagrange equation \cite{murray-rob}
\begin{equation*}
  \frac{\mathrm{d}}{\mathrm{d}t} \deriv[L]{\dq} - \deriv[L]{q} = F(q,\dq,u)
\end{equation*}
where $q$ is the system's generalized coordinates, $u$ represents the
external inputs (e.g, motor torque), $L$ is the Lagrangian (typically
kinetic energy minus potential energy for finite-dimensional
mechanical systems), and $F$ is the forcing function
that expresses external forces in the generalized coordinates.

\begin{figure}[h!]
  \centering
  \subfloat[Continuous Action Integral]{
    \label{fig:Action:Continuous}
    \begin{minipage}[b]{0.68\columnwidth}
      \centering\includegraphics[width=\textwidth]{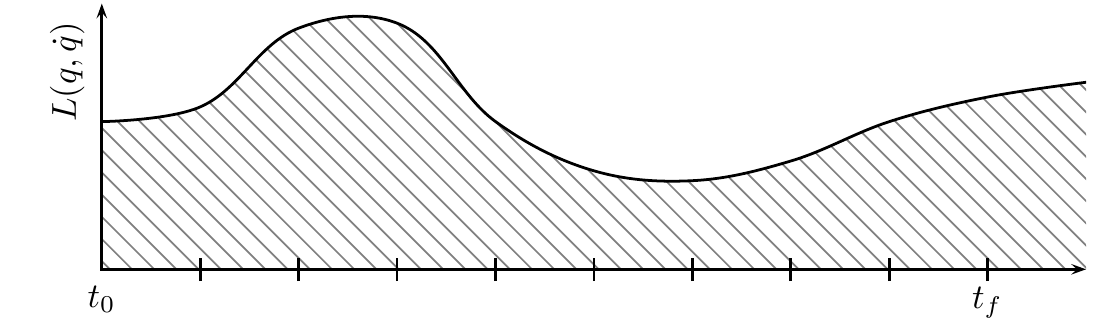}
    \end{minipage}}\\[12pt]
  \subfloat[Discrete Action Sum]{
    \label{fig:Action:Discrete}
    \begin{minipage}[b]{0.68\columnwidth}
      \centering\includegraphics[width=\textwidth]{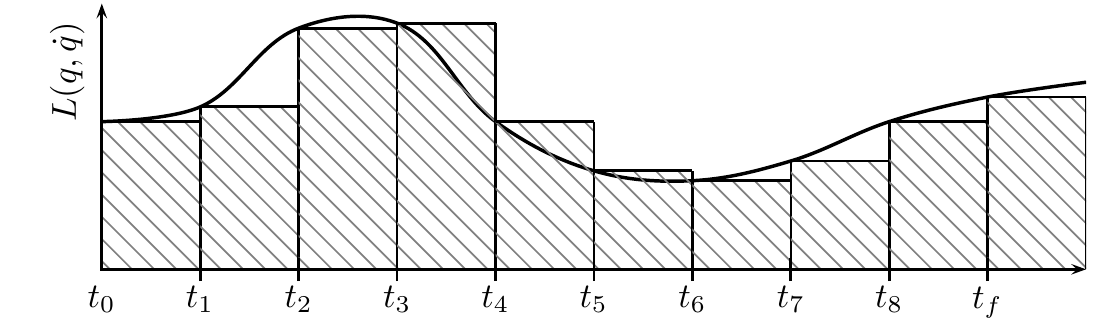}
  \end{minipage}}
  \caption{The continuous Euler-Lagrange equation is derived by
    minimizing the action integral \protect\subref{fig:Action:Continuous}.
    The discrete Euler-Lagrange equation is derived by minimizing the
    approximating action sum \protect\subref{fig:Action:Discrete}.}
 \label{fig:Action}
\vspace{-0.2in}
\end{figure}
The Euler-Lagrange equations can be derived from extremizing the
action integral, typically referred to as the (least) action principle.
The action integral---the integral of the Lagrangian with respect to
time along an arbitrary curve in the tangent bundle---is illustrated
as the shaded region in Fig.~\ref{fig:Action:Continuous}.  The action
principle stipulates that a mechanical system will follow the
trajectory that extremizes the action with respect to variations in
$q(t)$.  Applying calculus of variations to the action integral shows
that it is extremized by the Euler-Lagrange equation.

A variational integrator is derived by choosing a discrete Lagrangian,
$L_d$ that approximates the action over a discrete time step:
\begin{equation*}
  L_d(q_k, q_{k+1}) \approx \int_{t_k}^{t_{k+1}} L(q(s), \dq(s)) \ud s
\end{equation*}
where $q_k$ is a discrete-time configuration that approximates the
trajectory (i.e, $q_k \approx q(t_k)$).  This approximation can be achieved
with any quadrature rule; more accurate approximations lead to more
accurate integrators\cite{MarsdenWest2001}.  A concrete example of a
discrete Lagrangian approximation is in Sec.~\ref{sec:pend:DEL}.

By summing the discrete Lagrangian over an arbitrary trajectory, the
action integral is approximated by a discrete action, as shown in
Fig.~\ref{fig:Action}.  The action principle is then applied to
the action sum to find the {\it discrete trajectory} that extremizes
the discrete action.  The result of this calculation is the discrete
Euler-Lagrange (DEL) equations:
\begin{equation*}
  D_2L_d(q_{k-1}, q_k) + D_1L_d(q_k, q_{k+1}) = 0
\end{equation*}
where $D_nL_d$ is the {\it slot derivative}\footnote{The slot
  derivative $D_nL(A_1, A_2, \dots)$ represents the derivative of the
  function $L$ with respect to the $n$-th argument, $A_n$.  In many
  cases, the arguments to the function $L$ will be dropped for clarity
  and compactness.  Hence, it is helpful to keep in mind that the slot
  derivative applies to the argument order provided in a function's
  definition.} of $L_d$.

The DEL equation depends on the previous, current, and future
configuration (but it does not depend on the velocity, making this
integrator an appealing representation of dynamics for embedded
systems that measure configurations but not velocities).  The DEL
equation can also be written in an equivalent {\it position-momentum}
form that only depends on the current and future time steps:
\begin{subequations}
\begin{gather}
  p_k + D_1L_d(q_k, q_{k+1}) = 0
  \label{equ:del-position-momentum-1}
  \\
  p_{k+1} = D_2L_d(q_k, q_{k+1})
  \label{equ:del-position-momentum-2}
\end{gather}
\label{equ:del-position-momentum}
\end{subequations}
where $p_k$ is the discrete momentum of the system at time $k$.  (By
these definitions, it should be clear that $-D_1L_d(q_k, q_{k+1})$ and
$D_2L_d(q_k, q_{k+1})$ are both playing the role of a Legendre
transform in discrete time, and are accordingly referred to as the
left and right Legendre transforms.)

Equation~\eqref{equ:del-position-momentum} imposes a constraint on the
current and future positions and momenta.  Given an initial state
$p_k$ and $q_k$, the implicit equation
\eqref{equ:del-position-momentum-1} is solved numerically to find the
next configuration $q_{k+1}$.  In general,
\eqref{equ:del-position-momentum-1} is a non-linear equation that
cannot be solved explicitly for $q_{k+1}$.  In practice, the equation
is solved using a numeric method such as the Newton-Raphson algorithm.
The next momentum is then calculated explicitly by
\eqref{equ:del-position-momentum-2}.  After an update, $k$ is
incremented and the process is repeated to simulate the system for as
many time steps as desired.

Variational integrators can be extended to include non-conservative
forcing (e.g., a motor torque or damping) by using a discrete form of
the Lagrange-d'Alembert principle \cite{Ober2011}.  The continuous
force is approximated by a left and right discrete force, $F_d^-$ and
$F_d^{+}$:
\begin{align*}
\int_{t_k}^{t_{k+1}} &F(q(s),\dq(s),u(s)) \cdot \delta q \ud s \\ & \approx
  F_d^-(q_k, q_{k+1}, u_k) \cdot \delta q_k + 
  F_d^+(q_k, q_{k+1}, u_k) \cdot \delta q_{k+1}
\end{align*}
where $u_k$ is the discretization of the continuous force inputs: $u_k \approx
u(t_k)$.  As with the discrete Lagrangian, the discrete forcing can be
approximated by any quadrature rule. Certain quadrature rules may result in
$F_{d}^{\pm}$ also being a function of $u_{k+1}$ e.g. see \cite{Ober2011}.  As
we want to eventually apply classic control analysis and synthesis techniques,
we restrict our choice of quadrature rules to those where $F_{d}^{\pm}$ is
independent of $u_{k+1}$, largely to keep the resulting linearization of the
dynamics causal with respect to the input $u$. A concrete example is presented
in Sec.~\ref{sec:pend:DEL}.

For clarity, we use the following abbreviations for the discrete
Lagrangian and discrete forces throughout this paper:
\begin{align*}
  L_k &= L_d\left(q_{k-1}, q_k\right) \\
  F_k^\pm &= F_d^\pm\left(q_{k-1}, q_k, u_{k-1}\right).
\end{align*}

The forced DEL equations are then:
\begin{subequations}
  \label{equ:discrete-euler-lagrange-momentum}
  \begin{gather}
    \label{equ:discrete-euler-lagrange-momentum-1}
    p_k + D_1L_{k+1} + \Fm_{k+1} = 0
    \\
    \label{equ:discrete-euler-lagrange-momentum-2}
    p_{k+1} = D_2L_{k+1} + \Fp_{k+1}.
  \end{gather}
\end{subequations}
Again, \eqref{equ:discrete-euler-lagrange-momentum} provides a way to
calculate the configuration and momentum at the next time step from
the current time step.  Given the previous state ($p_k$ and $q_k$) and
the current input ($u_k$), the next configuration is found by
implicitly solving \eqref{equ:discrete-euler-lagrange-momentum-1}.
The momentum at the next time step is then calculated explicitly by
\eqref{equ:discrete-euler-lagrange-momentum-2}.

In the following section, we provide an example of a
variational integrator for a simple one dimensional system.  We will
use this example to help keep the calculations as concrete as possible
during the development of the structured linearization results.


\subsection{Example: Pendulum \label{sec:pend:DEL}}

\begin{figure}[h!]
  \centering
  \includegraphics[height=1.0in]{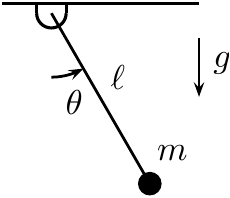}
  \caption{The pendulum is controlled by a torque at the pivot and subjected to
    the force of gravity $g$.}
  \label{fig:pendulum}
\end{figure}

Consider the pendulum shown in Fig.~\ref{fig:pendulum} with $m=\ell=1$ and a
gravitational force $g=9.8$. All units are assumed to be base units in SI.  The
pendulum has a single degree of freedom $\theta$, is controlled by a torque $u$
applied at the base, and is subjected to the acceleration due to gravity
$g$. The Lagrangian for the pendulum is
\begin{equation*}
  L(\theta, \dot{\theta}) = \half m \ell^2 \dot{\theta}^2 + m g \ell \cos\theta
                          = \half \dot{\theta}^2 + g \cos\theta.
\end{equation*}
The generalized force due to the torque input is:
\begin{equation*}
  F(\theta, \dot{\theta}, u) = u.
\end{equation*}

The discrete Lagrangian is found by approximating the
integral of the continuous-time Lagrangian over a short time interval
$\dt$ using the midpoint rule $\theta = \tfrac{\theta_k + \theta_{k+1}}{2}$ and
$\dot{\theta} = \tfrac{\theta_{k+1} - \theta_k}{\dt}$:

\begin{subequations}
  \begin{align}
    \label{equ:midpoint-dl}
    L_d(\theta_k, \theta_{k+1}) &= L\left(\tfrac{\theta_k + \theta_{k+1}}{2},
      \tfrac{\theta_{k+1} - \theta_k}{\dt}\right) \dt
    \\
    &= \tfrac{\left(\theta_{k+1} - \theta_k\right)^2}{2\dt} + g \dt \cos
    \tfrac{\theta_{k+1} + \theta_k}{2}.
  \end{align}
\end{subequations}
The forcing is approximated with a combination of a midpoint and
forward rectangle rule (though other choices of quadrature would be
fine as well):
\begin{align*}
  F^-_d(\theta_k, \theta_{k+1}, u_k) &= F(\tfrac{\theta_k + \theta_{k+1}}{2}, 
                                       \tfrac{\theta_{k+1} - \theta_k}{\dt}, u_k)\dt = u_k \dt
  \\
  F^+_d(\theta_k, \theta_{k+1}, u_k) &= 0.
\end{align*}

The first derivatives of $L_d$ are needed to implement the variational
integrator in \eqref{equ:discrete-euler-lagrange-momentum}:
\begin{align}
  D_1L_{k+1} &= - \tfrac{\theta_{k+1} - \theta_k}{\dt}
       - \tfrac{g \dt}{2} \sin \tfrac{\theta_{k+1} + \theta_k}{2} 
       \label{equ:ex2:D1Ld}
  \\
  D_2L_{k+1} &= \ \ \tfrac{\theta_{k+1} - \theta_k}{\dt}
       - \tfrac{g \dt}{2} \sin \tfrac{\theta_{k+1} + \theta_k}{2}.
       \label{equ:ex2:D2Ld}
\end{align}

The variational integrator update equations are found by
substituting \eqref{equ:ex2:D1Ld} into
\eqref{equ:discrete-euler-lagrange-momentum-1} and
\eqref{equ:ex2:D2Ld} into
\eqref{equ:discrete-euler-lagrange-momentum-2}:
\begin{subequations}
\begin{gather}
  p_k - \tfrac{\theta_{k+1} - \theta_k}{\dt}
       - \tfrac{g \dt}{2} \sin \tfrac{\theta_{k+1} + \theta_k}{2} + u_{k} \dt = 0
  \label{equ:pend:DEL-1}
  \\
  p_{k+1} = \tfrac{\theta_{k+1} - \theta_k}{\dt}
     - \tfrac{g \dt}{2} \sin \tfrac{\theta_{k+1} + \theta_k}{2}. 
  \label{equ:pend:DEL-2}
\end{gather}
\end{subequations}

Choose initial conditions $p_k = 0.5$, $q_k = \theta_k = 0.2$, a time step of
$\dt = 0.1s$, and an applied torque $u_k = 0.8$.  These values are substituted
in \eqref{equ:pend:DEL-1}, and a numeric root-finding algorithm finds the
unknown $\theta_{k+1}$.  In this case, the Newton-Raphson method was used to
find $\theta_{k+1} = 0.2471$.  Finally, the updated discrete momentum is
calculated using \eqref{equ:pend:DEL-2}: $p_{k+1} = 0.3627$. Computation of this
example utilizing {\tt trep} may be found at
\texttt{\url{https://trep.googlecode.com/git/examples/papers/tase2012/pend-single-step.py}}.
Note that it is already evident in this example that implicitly defined updates
are to be expected.  However, as we will see, these implicit updates have
explicit linearizations that can be computed as functions of the configuration.
Next we discuss the choice of state space for such a linearization.


\section{State Space Form\label{sec:state-form}}

In continuous time, the configuration and velocity of an
Euler-Lagrange system are often concatenated into a single state $x =
[q \;\; \dq ]^{T}$ to create a first-order
representation of the system.  This choice cannot be easily used for
the variational integrator because the finite-difference approximation
of the velocity involves configurations at different time steps.
Instead, the one-step representation of the integrator
\cite{GeometricNumericalIntegration} in
Eq.~(\ref{equ:discrete-euler-lagrange-momentum}) suggests that for the
variational integrator a convenient choice for the state is:
\begin{equation}
\label{eq-onestepupdate}
  x_{k+1} = \begin{bmatrix} q_{k+1} \\ p_{k+1} \end{bmatrix} = f(x_k,
  u_k), 
\end{equation}
where the function $f(x_k, u_k)$ is implicitly defined by
Eq.~(\ref{equ:discrete-euler-lagrange-momentum}).  However, the
Implicit Function Theorem guarantees that such a function exists
provided that the derivative
\begin{equation}
  M_{k+1} = D_2D_1L_{k+1} + D_2\Fm_{k+1}
  \label{equ:m-definition}
\end{equation}
is non-singular at $q_k$, $p_k$, and $u_k$.  This justifies
abstracting the discrete dynamics of the variational integrator this
way even though the underlying implementation still calculates the
update $q_{k+1}$ by numerically solving
\eqref{equ:discrete-euler-lagrange-momentum-1}.  The purpose of this
abstraction is to define the form for the linearization of the
discrete dynamics.  In the next section, we derive this linearization
and find that the derivatives of the abstract $f(x_k, u_k)$
representation are calculated explicitly.


\section{First Derivative \label{sec:first-deriv}}

Analysis and optimal control methods often rely on the first-order
linearization of system dynamics about a trajectory
\cite{AndersonMoore}.  The first-order linearization of the discrete
dynamics for the state form in Eq.~(\ref{eq-onestepupdate}) in
Sec.~\ref{sec:state-form} is:
\begin{gather}
  \notag
  \delta x_{k+1} = \deriv[f]{x_k} \delta x_k + 
     \deriv[f]{u_k} \delta u_k 
     \\
     \begin{bmatrix} \delta q_{k+1} \\ \delta p_{k+1} \end{bmatrix}
     = 
     \begin{bmatrix}
       \deriv[q_{k+1}]{q_k} & \deriv[q_{k+1}]{p_k} \\
       \deriv[p_{k+1}]{q_k} & \deriv[p_{k+1}]{p_k}
     \end{bmatrix}
     \begin{bmatrix} \delta q_k \\ \delta p_k \end{bmatrix}
     +
     \begin{bmatrix} \deriv[q_{k+1}]{u_k} \\ \deriv[p_{k+1}]{u_k} \end{bmatrix}
     \delta u_k.
     \label{equ:first-deriv}
\end{gather}

Six components are required to calculate this linearization.  These
derivatives are found directly from the variational integrator
equations \eqref{equ:discrete-euler-lagrange-momentum}, and all of
them result in explicit equations.

Derivatives of $q_{k+1}$ are found by implicitly differentiating
\eqref{equ:discrete-euler-lagrange-momentum-1} and solving for the
desired derivative.  We start by finding $\tderiv[q_{k+1}]{q_k}$:
\begin{align}
  &\tderiv{q_k} \left[ p_k + D_1L_{k+1} + \Fm_{k+1} = 0 \right]
  \notag
  \\
  &0 + D_1D_1L_{k+1} + D_2D_1L_{k+1} \tderiv[q_{k+1}]{q_k} 
   + D_1\Fm_{k+1}  \notag \\
  &\hspace{2.2in}  + D_2\Fm_{k+1} \tderiv[q_{k+1}]{q_k} = 0
  \notag
  \\
  &\left[D_2D_1L_{k+1} + D_2\Fm_{k+1}\right] \tderiv[q_{k+1}]{q_k} 
  = - \left[ D_1D_1L_{k+1} + D_1\Fm_{k+1}\right]
  \notag
  \\
  \label{equ:first-deriv-q-q}
  &\tderiv[q_{k+1}]{q_k}
  = - M_{k+1}^{-1} \left[ D_1D_1L_{k+1} + D_1\Fm_{k+1}\right]  
\end{align}
where $M_{k+1}$ is as defined by \eqref{equ:m-definition} and is
assumed to be non-singular (otherwise the Implicit Function Theorem
would not apply, making the state representation invalid).

The process is repeated to calculate $\tderiv[q_{k+1}]{p_k}$ and
$\tderiv[q_{k+1}]{u_k}$:
\begin{gather}
  \label{equ:first-deriv-q-p}
  \tderiv[q_{k+1}]{p_k} = -M_{k+1}^{-1}
  \\
  \label{equ:first-deriv-q-u}
  \tderiv[q_{k+1}]{u_k} = -M_{k+1}^{-1} \cdot D_3\Fm_{k+1}.
\end{gather}

Notice that each of these derivatives depends on the new configuration
$q_{k+1}$ (e.g, $D_1D_1L_{k+1} = D_1D_1L_d\left(q_k, q_{k+1}\right)$).
Before evaluating the derivatives, $q_{k+1}$ must be found by solving
\eqref{equ:discrete-euler-lagrange-momentum-1}.

Derivatives of $p_{k+1}$ are found directly by differentiating
\eqref{equ:discrete-euler-lagrange-momentum-2}:
\begin{gather}
  \tderiv[p_{k+1}]{q_k} = 
     \left[D_2D_2L_{k+1} + D_2\Fp_{k+1}\right]\tderiv[q_{k+1}]{q_k} +
     \notag \\
     D_1D_2L_{k+1} + D_1\Fp_{k+1} 
  \label{equ:first-deriv-p-q}
   \\
  \label{equ:first-deriv-p-p}
  \tderiv[p_{k+1}]{p_k} = \left[D_2D_2L_{k+1} + D_2\Fp_{k+1}\right] \tderiv[q_{k+1}]{p_k}
  \\
  \label{equ:first-deriv-p-u}
  \tderiv[p_{k+1}]{u_k} = \left[D_2D_2L_{k+1} + D_2\Fp_{k+1}\right] \tderiv[q_{k+1}]{u_k}
  + D_3\Fp_{k+1}.
\end{gather}

These derivatives depend on
\eqref{equ:first-deriv-q-q}--\eqref{equ:first-deriv-q-u}, so 
\eqref{equ:first-deriv-q-q}--\eqref{equ:first-deriv-q-u}
must be evaluated first.  Once calculated, their values are used in
\eqref{equ:first-deriv-p-q}--\eqref{equ:first-deriv-p-u} along with
the known value of $q_{k+1}$ to find the derivatives of $p_{k+1}$.
Once all six derivatives are calculated, they are organized into the
two matrices in \eqref{equ:first-deriv} to get the complete
first-order linearization about the current state.  Lastly, note that
the linearization is expressed entirely in terms of the discrete
Lagrangian's dependence on the configuration and the discrete forcing
function's dependence on the configuration and the continuous-time force.
This is critical in understanding how to calculate the linearization
without resorting to symbolic software, as discussed in
Section~\ref{sec:implementation}.


\subsection{Example: Pendulum {\it (cont.)}\label{sec:pend:first-deriv}}

We continue the pendulum example from \ref{sec:pend:DEL} by
calculating the first linearization (again, at the initial conditions
$p_k = 0.5$, $q_k = \theta_k = 0.2$, with a timestep of $\dt = 0.1$, and
an applied torque of $u_k = 0.8$).  The derivatives of the discrete
Lagrangian $L_d$ are:
\begin{align*}
  D_1D_1L_d &= \ \ \tfrac{1}{\dt}
       - \tfrac{g \dt}{4} \cos \tfrac{\theta_{k+1} + \theta_k}{2} = 9.7610
  \\
  D_2D_1L_d &= -\tfrac{1}{\dt}
       - \tfrac{g \dt}{4} \cos \tfrac{\theta_{k+1} + \theta_k}{2} = -10.2389
  \\
  D_1D_2L_d &= -\tfrac{1}{\dt}
       - \tfrac{g \dt}{4} \cos \tfrac{\theta_{k+1} + \theta_k}{2} = -10.2389
  \\
  D_2D_2L_d &= \ \ \tfrac{1}{\dt}
       - \tfrac{g \dt}{4} \cos \tfrac{\theta_{k+1} + \theta_k}{2} = 9.7610.
\end{align*}
The derivatives of the discrete forcing are trivial:
\begin{gather*}
  D_1F^-_d = D_2F^-_d = 0 
  \\
  D_3F^-_d = \dt.
\end{gather*}
Using these values with \eqref{equ:m-definition}, we find
$M_{k+1}^{-1} = (-10.2389 + 0)^{-1} = -0.0976$.  These are used with
\eqref{equ:first-deriv-q-q}--\eqref{equ:first-deriv-q-u} to calculate
the derivatives of $q_{k+1}$:
\begin{align*}
  \tderiv[q_{k+1}]{q_k} &=  0.0976 \cdot \left(9.7610 + 0\right) = 0.9533
  \\
  \tderiv[q_{k+1}]{p_k} &=  0.0976
  \\
  \tderiv[q_{k+1}]{u_k} &= 0.0976 \cdot 0.01 = 0.00976.
\end{align*}
These values are part of the linearization, but are also required to
calculate the derivatives of $p_{k+1}$ from
\eqref{equ:first-deriv-p-q}--\eqref{equ:first-deriv-p-u}.
\begin{align*}
  \tderiv[p_{k+1}]{q_k} &= (9.7610 + 0) \cdot 0.9533 + -10.2389 + 0 = -0.9333
  \\
  \tderiv[p_{k+1}]{p_k} &= (9.7610 + 0) \cdot 0.0976 = 0.9533
  \\
  \tderiv[p_{k+1}]{u_k} &= (9.7610 + 0) \cdot 0.00976 + 0 = 0.09533
\end{align*}
The six values define the entire first-order linearization:
\begin{equation*}
  \delta x_{k+1} = 
  \begin{bmatrix} 0.9533 & 0.0976 \\ -0.9333 & 0.9533 \end{bmatrix} \delta x_k
  + \begin{bmatrix} 0.00976 \\ 0.09533 \end{bmatrix} \delta u_k.
\end{equation*}

The first-order linearization frequently appears in analysis
applications.  For example, we can examine the controllability matrix of the
pendulum at this configuration to verify that it is linearly controllable:
\begin{gather*}
  \mathcal{C} = 
  \begin{bmatrix} B & AB \end{bmatrix} = 
  \begin{bmatrix} 0.00976 & 0.0186 \\ 0.09533 &  0.0818 \end{bmatrix}
  \\
  \mathrm{rank}\left(\mathcal{C}\right) = 2.
\end{gather*}
This linearization is carried out using {\tt trep} at
\texttt{\url{https://trep.googlecode.com/git/examples/papers/tase2012/pend-linearization.py}}.


\section{Second Derivative\label{sec:second-deriv}}

Optimal control applications can make use of the second-order linearization of
the dynamics to improve their convergence rate \cite{HauserSaccon}; this is
illustrated in Sec.~\ref{sec:example:pend:optimization}. The approach used in
Sec.~\ref{sec:first-deriv} extends to the second derivative of the dynamics as
well (we will call this the second-order linearization).  The expanded
second-order linearization of the discrete dynamics is
\begin{gather}
  \notag
  \delta^2 x_{k+1} = 
  \begin{bmatrix} \delta q_k \\ \delta p_k \\ \delta u_k \end{bmatrix}^T
  \begin{bmatrix}
    \derivII[f]{q_k}{q_k} & \derivII[f]{q_k}{p_k} & \derivII[f]{q_k}{u_k} \\[1em]
    \derivII[f^T]{q_k}{p_k} & \derivII[f]{p_k}{p_k} & \derivII[f]{p_k}{u_k} \\[1em]
    \derivII[f^T]{q_k}{u_k} & \derivII[f^T]{p_k}{u_k} & \derivII[f]{u_k}{u_k} \\[1em]
  \end{bmatrix}
  \begin{bmatrix} \delta q_k \\ \delta p_k \\ \delta u_k \end{bmatrix}
  \label{equ:second-deriv-definition}
\end{gather}
where symmetry is used to reduce the number of unique entries in the 
$3 \times 3$ array of third-order tensors
to six.  From Eq. \eqref{eq-onestepupdate} we have $f =
[q_{k+1} \;\; p_{k+1}]^{T}$, so each derivative of $f$ is calculated as two
third-order tensor components; one for $q_{k+1}$ and one for $p_{k+1}$. Thus 12
unique derivatives are needed for the second-order linearization.


\subsection{Derivatives of $q_{k+1}$}

As with the first derivatives, the second derivatives of $q_{k+1}$ are found by
differentiating \eqref{equ:discrete-euler-lagrange-momentum-1} twice and solving
for the desired derivative.  We will be using the notation $M\circ(X,Y)$ to
represent a bilinear operator\footnote{This is equivalent to the matrix
  representation $X^T M Y$ in finite dimensions, but this notation extends to
  infinite dimensional spaces such as those encountered in continuous trajectory
  optimization.} $M$ operating on $X$ and $Y$. We find
$\tderivII[q_{k+1}]{q_k}{q_k}$ as an example:
\begin{align*}
&  \tderivII{q_k}{q_k} \left[ p_k + D_1L_{k+1} + \Fm_{k+1} = 0 \right]
  \\
&  \tderiv{q_k} \left( 
    \left[D_2D_1L_{k+1} + D_2\Fm_{k+1}\right] \tderiv[q_{k+1}]{q_k} 
    =\right.\\
& \left. \hspace{1.5in}  - \left[ D_1D_1L_{k+1} + D_1\Fm_{k+1}\right]
    \right)
\end{align*} 
so we get that
\begin{align*}
&    \left[D_2D_1L_{k+1} + D_2\Fm_{k+1}\right]
    \tderivII[q_{k+1}]{q_k}{q_k} \\
& \hspace{0.5in}    + \left[D_1D_2D_1L_{k+1} + D_1D_2\Fm_{k+1}\right] \tderiv[q_{k+1}]{q_k} 
    \\
& \hspace{0.5in}   + \left[D_2D_2D_1L_{k+1} + D_2D_2\Fm_{k+1}\right] \op \left(\tderiv[q_{k+1}]{q_k}, \tderiv[q_{k+1}]{q_k} \right) 
    \\
&    = 
    -\left[ D_1D_1D_1L_{k+1} + D_1D_1\Fm_{k+1}\right] \\
& \hspace{0.5in}   -\left[ D_2D_1D_1L_{k+1} + D_2D_1\Fm_{k+1}\right] \tderiv[q_{k+1}]{q_k}.
\end{align*}
Solving for the desired derivative and substituting \eqref{equ:m-definition} we
get 
\begin{multline}
  \tderivII[q_{k+1}]{q_k}{q_k} = -M_{k+1}^{-1} \bigg(
    \left[ D_1D_1D_1L_{k+1} + D_1D_1\Fm_{k+1}\right]+
    \\
    \begin{aligned}
      \big[ D_1D_2D_1L_{k+1} + &D_1D_2\Fm_{k+1} + \\
      &D_2D_1D_1L_{k+1} + D_2D_1\Fm_{k+1} \big] \tderiv[q_{k+1}]{q_k}
    \end{aligned}
    \\
    + \left[D_2D_2D_1L_{k+1} + D_2D_2\Fm_{k+1}\right] \op \left(\tderiv[q_{k+1}]{q_k}, \tderiv[q_{k+1}]{q_k} \right) 
    \bigg)
  \label{equ:q-dq-dq}
\end{multline}
  
Previously we saw that the next state $x_{k+1}$ had to be found in order to
calculate the first derivatives.  Here we see that the second derivative
requires the first derivative as well.  This establishes the required order for
these calculations: The next state is found by the variational integrator, then
that state is used to calculate the first derivative, and then both results are
used to calculate the second derivative.

The other five second derivatives of $q_{k+1}$ are found by the same procedure.
The remaining derivatives with respect to state variables are:
\begin{multline}
  \tderivII[q_{k+1}]{q_k}{p_k}
  = 
  - M_{k+1}^{-1} \Big(
  \left[ D_2D_1D_1L_{k+1} + D_2D_1\Fm_{k+1}\right] \tderiv[q_{k+1}]{p_k}
  \\
  + \left[ D_2D_2D_1L_{k+1} + D_2D_2\Fm_{k+1} \right] \op \left(\tderiv[q_{k+1}]{q_k}, \tderiv[q_{k+1}]{p_k}\right) 
  \Big)
  \label{equ:q-dq-dp}
\end{multline}
and
\begin{multline}
  \tderivII[q_{k+1}]{p_k}{p_k} 
  = \\
  - M_{k+1}^{-1} \left[ D_2D_2D_1L_{k+1} + D_2D_2\Fm_{k+1} \right]
    \op \left( \tderiv[q_{k+1}]{p_k}, \tderiv[q_{k+1}]{p_k} \right).
  \label{equ:q-dp-dp}
\end{multline}

The derivatives with respect to state and input variables are:
\begin{multline}
  \tderivII[q_{k+1}]{q_k}{u_k}
  = 
  -M_{k+1}^{-1} \bigg(
  D_3D_1\Fm_{k+1}
  + D_3D_2\Fm_{k+1} \tderiv[q_{k+1}]{q_k}  
  \\
  + \Big[ D_2D_1D_1L_{k+1} + D_2D_1\Fm_{k+1} \Big] \tderiv[q_{k+1}]{u_k}
  \\
  + \left[ D_2D_2D_1L_{k+1} + D_2D_2\Fm_{k+1}\right] \op 
  \left(\tderiv[q_{k+1}]{q_k}, \tderiv[q_{k+1}]{u_k} \right)
  \bigg)
  \label{equ:q-dq-du}
\end{multline}
and
\begin{multline}
  \tderivII[q_{k+1}]{p_k}{u_k} 
  = 
  - M_{k+1}^{-1} \bigg(
  D_3D_2\Fm_{k+1} \tderiv[q_{k+1}]{p_k}  \\
  +\left[ D_2D_2D_1L_{k+1} + D_2D_2\Fm_{k+1} \right] \op 
  \left(\tderiv[q_{k+1}]{p_k}, \tderiv[q_{k+1}]{u_k} \right)
  \bigg).
  \label{equ:q-dp-du}
\end{multline}

The second derivative of the next configuration with respect to the
inputs is:
\begin{multline}
  \tderivII[q_{k+1}]{u_k}{u_k}
  =\\ 
  - M_{k+1}^{-1} \bigg(
  D_3D_3\Fm_{k+1} 
  + \left[ D_3D_2\Fm_{k+1} + D_2D_3\Fm_{k+1} \right] \tderiv[q_{k+1}]{u_k}
  \\
  + \left[ D_2D_2D_1L_{k+1} + D_2D_2\Fm_{k+1} \right]
    \op \left( \tderiv[q_{k+1}]{u_k}, \tderiv[q_{k+1}]{u_k} \right)
    \bigg).
  \label{equ:q-du-du}
\end{multline}

These six derivatives make up the entire second linearization of the
first half of the state.


\subsubsection{Pendulum {\it (cont.)}}
We now return to the pendulum example, again with the same initial data.
We calculate the first part of the second-order linearization directly
from~\eqref{equ:q-dq-dq}--\eqref{equ:q-du-du}.  The next state and
first-order linearization of the pendulum were already found in
Sec.~\ref{sec:pend:DEL} and \ref{sec:pend:first-deriv}, respectively.
After calculating the higher derivatives of $L_{k+1}$ and
$F^\pm_{k+1}$, we find the second derivative of the pendulum's
configuration dynamics $\delta^2 q_{k+1} = $
\begin{gather}
  \notag
  \begin{bmatrix} \delta q_k \\ \delta p_k \\ \delta u_k \end{bmatrix}^T
  \begin{bmatrix}
    1.01\e{-2} & 5.06\e{-4} & 5.06\e{-5} \\
    5.06\e{-4} & 2.53\e{-5} & 2.53\e{-6} \\
    5.06\e{-5} & 2.53\e{-6} & 2.53\e{-7} 
  \end{bmatrix}
  \begin{bmatrix} \delta q_k \\ \delta p_k \\ \delta u_k \end{bmatrix}.
\end{gather}

We also need the second derivatives of $p_{k+1}$ to complete the
second-order linearization.


\subsection{Derivatives of $p_{k+1}$}

The six derivatives of the momentum component of the state, $p_{k+1}$
are calculated directly from the explicit momentum equation
\eqref{equ:discrete-euler-lagrange-momentum-2}.  The derivatives with
respect to state variables are:
\begin{multline}
  \tderivII[p_{k+1}]{q_k}{q_k} = 
  D_1D_1D_2L_{k+1} 
  + D_1D_1\Fp_{k+1} \\
  + \left[ D_2D_2L_{k+1} + D_2\Fp_{k+1} \right] \tderivII[q_{k+1}]{q_k}{q_k}  
  \\
  \begin{aligned}
    + \big[ D_2D_1D_2L_{k+1} &+ D_1D_2D_2L_{k+1} + \\ 
    D_2D_1\Fp_{k+1} &+ D_1D_2\Fp_{k+1} \big] \tderiv[q_{k+1}]{q_k}
  \end{aligned}
  \\
  + \left[D_2D_2D_2L_{k+1} + D_2D_2\Fp_{k+1} \right] 
    \op \left( \tderiv[q_{k+1}]{q_k}, \tderiv[q_{k+1}]{q_k} \right)
    \label{equ:p-dq-dq}
\end{multline}
\begin{multline}
  \tderivII[p_{k+1}]{q_k}{p_k}
  = 
  \left[ D_2D_1D_2L_{k+1} + D_2D_1\Fp_{k+1}\right] \tderiv[q_{k+1}]{p_k}
  \\
  + \big[ D_2D_2L_{k+1} + D_2\Fp_{k+1} \big] \tderivII[q_{k+1}]{q_k}{p_k}
  \\
  + \left[ D_2D_2D_2L_{k+1} + D_2D_2\Fp_{k+1} \right] 
    \op \left( \tderiv[q_{k+1}]{q_k}, \tderiv[q_{k+1}]{p_k}\right)
    \label{equ:p-dq-dp}
\end{multline}
\begin{multline}
  \tderivII[p_{k+1}]{p_k}{p_k} = 
  \left[ D_2D_2L_{k+1} + D_2\Fp_{k+1} \right] \tderivII[q_{k+1}]{p_k}{p_k}
  \\
  + \left[ D_2D_2D_2L_{k+1} + D_2D_2\Fp_{k+1}\right]
  \op \left(\tderiv[q_{k+1}]{p_k}, \tderiv[q_{k+1}]{p_k}\right).
    \label{equ:p-dp-dp}
\end{multline}

The derivatives with respect to state and input variables are:
\begin{multline}
  \tderivII[p_{k+1}]{q_k}{u_k}
  = 
  D_3D_1\Fp_{k+1} 
  + D_3D_2\Fp_{k+1} \tderiv[q_{k+1}]{q_k}
  + \big[ D_2D_1D_2L_{k+1} \\ + D_2D_1\Fp_{k+1} \big] \tderiv[q_{k+1}]{u_k} 
  + \left[ D_2D_2L_{k+1} + D_2\Fp_{k+1} \right] \tderivII[q_{k+1}]{q_k}{u_k}
  \\
  + \left[ D_2D_2D_2L_{k+1} + D_2D_2\Fp_{k+1} \right]
    \op \left( \tderiv[q_{k+1}]{q_k}, \tderiv[q_{k+1}]{u_k} \right)
    \label{equ:p-dq-du}
\end{multline}
\begin{multline}
  \tderivII[p_{k+1}]{p_k}{u_k} = 
  D_3D_2\Fp_{k+1} \tderiv[q_{k+1}]{p_k} + \\
  \left[ D_2D_2L_{k+1} + D_2\Fp_{k+1} \right] \tderivII[q_{k+1}]{p_k}{u_k}
  \\
  + \left[ D_2D_2D_2L_{k+1} + D_2D_2\Fp_{k+1} \right] 
    \op \left( \tderiv[q_{k+1}]{p_k}, \tderiv[q_{k+1}]{u_k} \right).
    \label{equ:p-dp-du}
\end{multline}
Finally, the second derivative with respect to the input variables is:
\begin{multline}
  \tderivII[p_{k+1}]{u_k}{u_k}
  = 
  D_3D_3\Fp_{k+1}
  + \left[ D_3D_2\Fp_{k+1} + D_2D_3\Fp_{k+1} \right] \tderiv[q_{k+1}]{u_k}
  \\
  + \left[ D_2D_2D_2L_{k+1} + D_2D_2\Fp_{k+1} \right] 
    \op \left( \tderiv[q_{k+1}]{u_k}, \tderiv[q_{k+1}]{u_k} \right)
  \\
  + \big[ D_2D_2L_{k+1} + D_2\Fp_{k+1} \big] \tderivII[q_{k+1}]{u_k}{u_k}.
    \label{equ:p-du-du}
\end{multline}

As with the first derivatives, the second derivatives of $p_{k+1}$
depend on those of $q_{k+1}$.  We handle this dependency by first
evaluating \eqref{equ:q-dq-dq}--\eqref{equ:q-du-du} to get their
numerical values and then plug those values into
\eqref{equ:p-dq-dq}--\eqref{equ:p-du-du}.

\textbf{Note:} By evaluating each of the twelve equations above, we
explicitly calculate the complete second-order linearization for a
forced system in generalized coordinates.
Section~\ref{sec:constrained-systems} describes how this approach is
extended to systems with holonomic constraints.


\subsubsection{Pendulum {\it (cont.)}}

We complete the second derivative by evaluating
\eqref{equ:p-dq-dq}--\eqref{equ:p-du-du} with the values found earlier
and the remaining derivatives of $L_{k+1}$ and $F^\pm_{k+1}$.  The
result is $  \delta^2 p_{k+1} = $
\begin{gather}
  \notag
  \begin{bmatrix} \delta q_k \\ \delta p_k \\ \delta u_k \end{bmatrix}^T
  \begin{bmatrix}
    2.02\e{-1} & 1.01\e{-2} & 1.01\e{-3} \\
    1.01\e{-2} & 5.06\e{-4} & 5.06\e{-5} \\
    1.01\e{-3} & 5.06\e{-5} & 5.06\e{-6} 
  \end{bmatrix}
  \begin{bmatrix} \delta q_k \\ \delta p_k \\ \delta u_k \end{bmatrix}.
\end{gather}

These second-order linearizations $\delta^{2}q_{k+1}$ and $\delta^{2}p_{k+1}$
carried out using {\tt trep} can be found at
\texttt{\url{https://trep.googlecode.com/git/examples/papers/tase2012/pend-linearization.py}}.


\section{Constrained Systems\label{sec:constrained-systems}}

In this section, we discuss how the approach described in
Sec.~\ref{sec:first-deriv} and \ref{sec:second-deriv} to calculate the
discrete linearizations extends to constrained variational
integrators. Variational integrators are particularly well-suited to
systems with holonomic constraints because the update equation for a
constrained variational integrator explicitly incorporates the
holonomic constraint.  This is opposed to replacing the holonomic
constraint with its derivatives and then
projecting the update onto the feasible set, a common approach in
numeric integration of ordinary differential equations.
Variational integrators enforce the holonomic constraint at every time
step while still preserving the symplectic form and conserving
momentum.

Variational integrators for constrained systems\cite{MarsdenWest2001}
are derived using the same Lagrange-multiplier method used in the
continuous case\cite{murray-rob}.  Given a continuous-time constraint
of the form $h(q)=0$, the DEL equations for a forced, constrained
variational integrator are:
\begin{subequations}
  \label{equ:discrete-euler-lagrange-constrained}
  \begin{gather}
    \label{equ:discrete-euler-lagrange-constrained-1}
    p_k + D_1L_{k+1} + \Fm_{k+1} - Dh^T(q_k)\lambda_k = 0
    \\
    \label{equ:discrete-euler-lagrange-constrained-2}
    h(q_{k+1}) = 0
    \\
    \label{equ:discrete-euler-lagrange-constrained-3}
    p_{k+1} = D_2L_{k+1} + \Fp_{k+1}
  \end{gather}
\end{subequations}
where $\lambda_k$ are the Lagrange multipliers that can be interpreted as
discrete-time forces enforcing the constraint.  In this case, given
$p_k$ and $q_k$, a root-finding algorithm solves
\eqref{equ:discrete-euler-lagrange-constrained-1} and
\eqref{equ:discrete-euler-lagrange-constrained-2} to find $q_{k+1}$
and $\lambda_k$.  The updated momentum $p_{k+1}$ is then explicitly
calculated from \eqref{equ:discrete-euler-lagrange-constrained-3}.

The Lagrange multipliers are completely determined by $q_k$, $p_k$,
and $u_k$, so the state representation from
Section~\ref{sec:state-form} is unchanged.  Accordingly, the same
derivatives $\delta x_{k+1}$ and $\delta^2x_{k+1}$ are needed to find the
linearizations.  Rather than derive every equation, we calculate one
component of the first and second derivatives to demonstrate the
process.

For the first-order linearization, we find
$\tderiv[q_{k+1}]{q_k}$.  We start by differentiating
\eqref{equ:discrete-euler-lagrange-constrained-1}:
\begin{gather}
  \tderiv{q_k}\left[p_k + D_1L_{k+1} + \Fm_{k+1} - Dh^T(q_k)\lambda_k = 0\right]
  \notag
  \\
\Rightarrow \ \ \ \   \tderiv[q_{k+1}]{q_k}
  = - M_k^{-1} \Big[ C_{q_k}
    - Dh^T(q_k)\tderiv[\lambda_k]{q_k}
    \Big]
  \label{equ:constrained-first-deriv-q-q}
\end{gather}
where
\begin{equation*}
  C_{q_k} = 
    D_1D_1L_{k+1} 
    \\
    + D_1\Fm_{k+1} 
    - D^2h^T(q_k)\lambda_k.
\end{equation*}
To evaluate this derivative, we must calculate
$\tderiv[\lambda_k]{q_k}$.  This is found by differentiating
\eqref{equ:discrete-euler-lagrange-constrained-2}, substituting in
\eqref{equ:constrained-first-deriv-q-q}, and solving for
$\tderiv[\lambda_k]{q_k}$:
\begin{gather}
  \tderiv{q_k} \left[ h(q_{k+1}) = 0 \right]
  \notag
  \\
  Dh(q_{k+1}) \tderiv[q_{k+1}]{q_k} = 0
  \notag
  \\
  Dh(q_{k+1}) M_k^{-1} \Big[ C_{q_k}
    - Dh^T(q_k)\tderiv[\lambda_k]{q_k}
    \Big] = 0 
  \notag
  \\
  Dh(q_{k+1}) M_k^{-1} C_{q_k}
  - Dh(q_{k+1}) M_k^{-1} Dh^T(q_k)\tderiv[\lambda_k]{q_k} 
  = 0
  \notag
  \\
  \tderiv[\lambda_k]{q_k} 
  =
  \left[ Dh(q_{k+1}) M_k^{-1} Dh^T(q_k) \right]^{-1} 
  Dh(q_{k+1}) M_k^{-1} C_{q_k}.
  \label{equ:constrained-first-deriv-lambda-q}
\end{gather}  
To calculate $\tderiv[q_{k+1}]{q_k}$, the  constrained DEL equation
\eqref{equ:discrete-euler-lagrange-constrained} is solved numerically
to find $q_{k+1}$ and $\lambda_k$.  These values are used in
\eqref{equ:constrained-first-deriv-lambda-q} to find
$\tderiv[\lambda_k]{q_k}$.  Finally, $\tderiv[q_{k+1}]{q_k}$ is
calculated with \eqref{equ:constrained-first-deriv-q-q}.

The same approach is used to find the remaining components of the first
derivative, so we do not repeat the derivation here.  Derivations of the
  remaining components can be found in \cite{JohnsonThesis}. We continue onto
the second derivative by calculating $\tderivII[q_{k+1}]{q_k}{q_k}$.
\begin{gather}
  \tderivII{q_k}{q_k}\left[p_k + D_1L_{k+1} + \Fm_{k+1} - Dh^T(q_k)\lambda_k = 0\right]
  \notag
  \\
\Rightarrow\ \ \ \ \   \tderivII[q_{k+1}]{q_k}{q_k} = 
  - M_{k+1}^{-1} \bigg(   C_{q_kq_k} 
    - Dh^T(q_k)\tderivII[\lambda_k]{q_k}{q_k}
    \bigg)
    \label{equ:constrained-second-deriv-q-q-q}
\end{gather}
where
\begin{align*}
  C_{q_kq_k} =& 
      D_1D_1D_1L_{k+1}
    + D_1D_1\Fm_{k+1} 
    \\
    &+ \Big[ 
      D_2D_1D_1L_{k+1} 
      + D_2D_1\Fm_{k+1} \Big. \\ \Big.
       & \hspace{0.5in} + D_1D_2D_1L_{k+1}
       + D_1D_2\Fm_{k+1} 
       \Big] \tderiv[q_{k+1}]{q_k}
    \\
    &+ \Big[ D_2D_2D_1L_{k+1} + D_2D_2\Fm_{k+1} \Big] 
    \op \left( \tderiv[q_{k+1}]{q_k}, \tderiv[q_{k+1}]{q_k} \right)
    \\
    &- D^3h^T(q_k)\lambda_k 
    - 2 D^2h^T(q_k)\tderiv[\lambda_k]{q_k}.
\end{align*}
Again, we find the corresponding second derivative of $\lambda_k$ by
differentiating \eqref{equ:discrete-euler-lagrange-constrained-2}
twice:
\begin{gather*}
  \tderivII{q_k}{q_k} \left[ h(q_{k+1}) = 0 \right]
  \\
  D^2h(q_{k+1}) \op \left( \tderiv[q_{k+1}]{q_k}, \tderiv[q_{k+1}]{q_k} \right)
  + Dh(q_{k+1}) \tderivII[q_{k+1}]{q_k}{q_k} = 0.
\end{gather*}
We substitute in \eqref{equ:constrained-second-deriv-q-q-q} and solve for
$\tderivII[\lambda_k]{q_k}{q_k}$:
\begin{multline}
  \tderivII[\lambda_k]{q_k}{q_k}
  =
  \left[Dh(q_{k+1}) M_{k+1}^{-1} Dh^T(q_k) \right] \cdot
  \\
  \left[
    Dh(q_{k+1}) M_{k+1}^{-1} C_{q_kq_k}
    -  D^2h(q_{k+1}) \op \left( \tderiv[q_{k+1}]{q_k}, \tderiv[q_{k+1}]{q_k} \right)
  \right].
  \label{equ:constrained-second-deriv-lambda-q-q}
\end{multline}

To calculate this $\tderivII[q_{k+1}]{q_k}{q_k}$, we solve for the
next state, calculate the first derivatives, evaluate
\eqref{equ:constrained-second-deriv-lambda-q-q} to find
$\tderivII[\lambda_k]{q_k}{q_k}$, and finally evaluate
\eqref{equ:constrained-second-deriv-q-q-q} to find the second
derivative.  This same procedure is used to calculate the other
components of the constrained second derivative.

Note that the constrained momentum update
\eqref{equ:discrete-euler-lagrange-constrained-3} is identical to the
unconstrained case \eqref{equ:discrete-euler-lagrange-momentum-2}, so
the first- and second-order linearizations are identical.


\section{Singularities of the Linearization}
\label{sec:sing-linearization}
Eq. \eqref{equ:m-definition} includes a matrix $M_{k+1}$ that must be inverted
in both the first and second order linearizations presented in
Sec. \ref{sec:first-deriv} and Sec. \ref{sec:second-deriv} respectively.
Additionally Sec. \ref{sec:constrained-systems} shows that when the system
involves constraints, the linearizations additionally involve the term $\left[
  Dh(q_{k+1}) M_k^{-1} Dh^T(q_k) \right]^{-1}$.  For a general mechanical
system, the requirements for invertibility of these two terms are not known, but
certainly the choice of coordinate chart can cause singularities as can
degeneracy of the Lagrangian system.  We illustrate this point by presenting two
simple examples that demonstrate situations where $M_{k+1}$ and $\left[
  Dh(q_{k+1}) M_k^{-1} Dh^T(q_k) \right]$ become singular.

\subsection{Singularities of $M_{k+1}$ in a Spherical Pendulum}
\label{sec:sing-spherical-pend}
Consider an unforced, spherical pendulum of mass $m$ and length $r$
  under the influence of gravity $g$ in generalized coordinates $q = (\theta,
  \phi)$ where $\theta$ is the polar angle measured from the zenith direction
  which is aligned with gravity, and $\phi$ is the azimuthal angle. For this
  system, the Lagrangian is given by
\begin{equation*}
  L(q,\dot{q}) = \half m r^{2}\left(\dot{\theta}^{2} + 
    \sin^{2}{\theta} \dot{\phi}^{2}\right) + m g r \cos{\theta}.
\end{equation*}
It is well-known that this choice of generalized coordinates does not provide a
global chart, resulting in singular configurations. This can be seen by looking
at the mass matrix for this system which is given by
\begin{equation}
  \label{equ:sphere-pend-m}
  M(q) = \derivII[L]{\dot{q}}{\dot{q}} =
  \begin{bmatrix}
    mr^{2}\sin^{2}(\theta) & 0 \\
    0 & mr^{2}
  \end{bmatrix}.
\end{equation}
Clearly this matrix is singular if $\theta=n\pi \;\forall \;n \in \mathbb{Z}$
where $\mathbb{Z}$ is the set of all integers. Using the midpoint-rule discrete
Lagrangian of \eqref{equ:midpoint-dl} we can construct the discrete Lagrangian
for this system as
\begin{multline*}
  L_{d}(q_{k},q_{k+1}) = \tfrac{m r^{2}}{2} \left(\tfrac{\theta_{k+1} -
      \theta_{k}}{\Delta t}\right)^{2} + \\
  \tfrac{m r^{2}}{2} \sin^{2}\left(\tfrac{\theta_{k+1} + \theta_{k}}{2}\right)
  \left(\tfrac{\phi_{k+1} - \phi_{k}}{\Delta t}\right)^{2}
  + m g r \cos\left(\tfrac{\theta_{k+1} + \theta_{k}}{2}\right).
\end{multline*}
From \eqref{equ:m-definition} $M_{k+1}$ may be calculated as 
\begin{equation*}
  M_{k+1} =
  \begin{bmatrix}
    \tderivII[L_{k+1}]{\phi_{k+1}}{\phi_{k}} & \tderivII[L_{k+1}]{\theta_{k+1}}{\phi_{k}} \\[12pt]
    \tderivII[L_{k+1}]{\phi_{k+1}}{\theta_{k}} & \tderivII[L_{k+1}]{\theta_{k+1}}{\theta_{k}}
  \end{bmatrix}.
  \label{equ:mk1-spherical}
\end{equation*}
To find singularities of this matrix, set its determinant equal to zero.  The
determinant of $M_{k+1}$ is 
\begin{multline}
  \label{equ:spherical-m-det}
  \det(M_{k+1}) = \tfrac{m^{2}r^{3}}{4 \Delta t^{2}}
  \sin^{2}\left(\tfrac{\theta_{k+1} + \theta_{k}}{2}\right) \bigg[
  \Delta t^{2} g \cos \left(\tfrac{\theta_{k+1} + \theta_{k}}{2}\right) \\ +
  r\big(2 + \cos\left(\theta_{k+1} + \theta_{k}\right)\big) \left(\phi_{k} -
    \phi_{k+1}\right)^{2} + 4r\bigg].
\end{multline}

If the term before the brackets is zero then $\det(M_{k+1}) = 0$ and the matrix
is singular.  This implies $M_{k+1}$ is singular if $\theta_{k+1} + \theta_{k} =
n\pi \;\forall\; n \in \mathbb{Z}$.  Thus if the discrete system is at the
singular configuration of the continuous system for two consecutive timesteps or
if the consecutive configurations are symmetric about the singular
configuration, then $M_{k+1}$ is non-invertible.

Equation \eqref{equ:spherical-m-det} is also zero if the term in the brackets is
zero, implying that there exist sets of consecutive configurations that 
cause $M_{k+1}$ to be singular that are not directly related to the continuous-time
singular configurations. However, there  exists an
upper bound on the timestep size $\Delta t$ that prevents the bracketed term
from being zero.  To illustrate,  assume that the
spherical pendulum is in pure pendular motion i.e. $\phi_{k+1} = \phi_{k}$.
With this assumption setting the bracketed term to zero yields
\begin{align*} &\Delta t^{2} g \cos
    \left(\tfrac{\theta_{k+1} + \theta_{k}}{2}\right) + 4r =
    0 \\
    \implies &\theta_{k+1} + \theta_{k} = 2 \cos^{-1}\left(\frac{-4r}{\Delta
        t^{2}g}\right).
\end{align*}
The inverse cosine is only defined if its argument $\tfrac{-4r}{\Delta t^{2}g}
\in [-1,1]$.  Noting that $g$, $r$, and $\Delta t$ are all greater than zero, we
see that the argument is bounded above by zero.  Thus the inverse cosine only
has a solution if 
\begin{align*}
  &-1  \leq \frac{-4r}{\Delta t^{2}g} \\
  &\implies \Delta t \geq \sqrt{\frac{4r}{g}}.
\end{align*}
Thus if $\Delta t < \sqrt{(4r)/{g}}$ and $\phi_{k+1} = \phi_{k}$ the
bracketed term is always nonzero and the only singularity that exists in
$M_{k+1}$ is the one induced by our choice of generalized coordinates.

If we relax the constraint $\phi_{k+1} = \phi_{k}$ similar reasoning shows that
for a given set of constants $g$, $r$, and $\Delta t$ then there exists an upper
bound, $\lambda^{*}$, on the difference in $\phi$ across a timestep such that if
$\abs{\phi_{k} - \phi_{k+1}} \leq \lambda^{*}$ then the bracketed term cannot be
zero. Thus for a given spherical pendulum, as long as the timestep is small
enough, the only singularities that exist in $M_{k+1}$
are singularities that are caused by the choice of generalized coordinates.

This example illustrates that singularities caused by a poor choice
of local, generalized coordinates can show up as singularities in $M_{k+1}$.
Additionally, other singularities in $M_{k+1}$ can be related to
coarse sampling of the continuous system and can be prevented by decreasing the
mesh size.

\subsection{Singularities of a Constrained Pendulum}
\label{sec:sing-constr-pend}

Consider again a simple pendulum system as shown in Fig. \ref{fig:pendulum}.
Assume that gravity is zero.  In this example, we will represent the system in
Cartesian coordinates $q=(x,y)$, and we will add a constraint of the form
\begin{equation*}
  h(q) = x^{2} + y^{2} - l^{2} = 0.
\end{equation*}
This system's midpoint discrete Lagrangian is

\begin{equation*}
  L_{d}(q_{k}, q_{k+1}) = \Half m \left( \left(\tfrac{x_{k+1}-x_{k}}{\Delta
        t}\right)^{2} + \left(\tfrac{y_{k+1}-y_{k}}{\Delta t}\right)^{2} \right) .
\end{equation*}
We calculate $M_{k+1}$ as
\begin{equation*}
  M_{k+1} =
  \begin{bmatrix}
    \tderivII[L_{k+1}]{x_{k+1}}{x_{k}} & \tderivII[L_{k+1}]{y_{k+1}}{x_{k}} \\[12pt]
    \tderivII[L_{k+1}]{x_{k+1}}{y_{k}} & \tderivII[L_{k+1}]{y_{k+1}}{y_{k}}
  \end{bmatrix} = 
  \begin{bmatrix}
    \tfrac{-m}{\Delta t^{2}} & 0 \\[12pt]
    0 & \tfrac{-m}{\Delta t^{2}}
  \end{bmatrix}
  \label{equ:mk1-constrained-pend}
\end{equation*}
and see that it is always invertible.  In order to linearize this system $\left[
  Dh(q_{k+1}) M_k^{-1} Dh^T(q_k) \right]$ must be invertible. Using the given
$h(q)$
\begin{equation*}
  \left[Dh(q_{k+1}) M_k^{-1} Dh^T(q_k) \right] = -\frac{4\Delta t^{2}}{m}(x_{k+1}x_{k} +
  y_{k+1}y_{k}).
\end{equation*}

This is only non-invertible if $Dh(q_{k})$ and $Dh(q_{k+1})$ are orthogonal in
the Euclidean sense.  This orthogonality would require the pendulum to move
$\pi/2~rad$ in a single timestep.  Regardless of how fast the pendulum is
moving, we are guaranteed that as $\Delta t \rightarrow 0$ the change in
configuration goes to zero i.e. $\|q_{k} - q_{k+1}\| \rightarrow 0$.  Thus there
is always a sufficiently small timestep to ensure $\left[Dh(q_{k+1}) M_k^{-1}
  Dh^T(q_k) \right]$ is non-singular. From this example, it is clear that for a
general mechanical system, even in cases where $M_{k}$ is non-singular there may
be situations where the term $\left[Dh(q_{k+1}) M_k^{-1} Dh^T(q_k) \right]$ may
be singular.

Generalizing the ways in which a linearization of a variational integrator can
cease to be well-posed is clearly an important issue to pursue.  However, for
purposes of a numerical method, it suffices to check that $M_k$ and
$\left[Dh(q_{k+1}) M_k^{-1} Dh^T(q_k) \right]$ are invertible at every time step
to be confident that the computed linearization is correct.


\section{Implementation\label{sec:implementation}}

Deriving the first- and second-order derivative equations of
Sec. \ref{sec:first-deriv} and \ref{sec:second-deriv} is  procedurally
straightforward.  The more complicated issue in implementing the  first- and second-order derivatives is
calculating higher derivatives of the discrete Lagrangian and discrete forces.
Approaches relying on the symbolic equations, as was done for the pendulum
example, do not scale to complex systems even with the help of symbolic algebra
software.

In \cite{JohnsonMurpheyTRO2008}, the authors describe a method to implement
variational integrators using a hierarchical tree representation
\cite{featherstone,NakamuraYamane,YamaneBook}.  That approach calculates exact
derivatives of the discrete Lagrangian numerically and scales to large, complex
mechanical systems like a biomechanical model of the human
hand\cite{JohnsonWAFR2010}.  Moreover, it calculates the derivatives of the
discrete Lagrangian without ever having an explicit representation of the
discrete Lagrangian itself; all the calculations are implicitly defined in terms
of the tree representation that encodes the mechanical topology.  

It should be noted that \emph{how} one uses the mechanical topology to compute
derivatives is not critical for the present work.  Indeed, one could use any
method that allows one to compute higher-order derivatives of the discrete
Lagrangian (e.g., a spatial-operator approach \cite{Jain} or a recursive
dynamics \cite{Anderson} approach).  So long as the computational method
provides the derivatives from Section~\ref{sec:first-deriv} and
\ref{sec:second-deriv}, the only risk is that one method may provide better or
worse scaling properties than another method, something that is beyond the scope
of the present paper.

The method in \cite{JohnsonMurpheyTRO2008} naturally extends to calculating
higher-order derivatives\cite{JohnsonACC2010} of the discrete Lagrangian and
discrete forcing that are needed by the first- and second-order linearizations
\cite{MurpheyACC2011}.  Thus it provides a complete framework for implementing
variational integrators and their first- and second-order linearizations.

These algorithms have been implemented in an open source software library called
{\tt trep}.  The software calculates continuous and discrete dynamics, along
with their first- and second-order linearizations, for arbitrary mechanical
systems in generalized coordinates, including those with holonomic constraints.
{\tt trep} is freely available at \texttt{\url{http://trep.googlecode.com}}.

The following examples demonstrate the application of the first-order
linearization as a tool for generating stabilizing controllers; the first
example considers a simple pendulum, and the second features a more complex
mechanical system highlighting the scalability of the method described in
\cite{JohnsonMurpheyTRO2008}.  Additionally, the first example utilizes a
second-order optimization technique requiring the second-order linearization
presented in Sec. \ref{sec:second-deriv}. All of the calculations were performed
by {\tt trep}.

\subsection{Example: Optimal Control of the Pendulum}
\label{sec:example:pend:optimization}
As a simple example, we once again consider the pendulum in
Fig. \ref{fig:pendulum}.  Using the choice of state discussed in
Sec. \ref{sec:state-form} the state for this system is $x(k) = [q(k) \;\;
p(k)]^{T}$.  Here we are dropping the shortcut of using a subscript $k$ to
indicate sequence index to avoid confusion with other subscripts used for
distinguishing trajectories.  We begin by defining a dynamically infeasible,
discrete reference trajectory over the time horizon $t_{ref} = \{t_{ref}(k) = k
\Delta t \;|\; k=0,\ldots,N\}$ with $\Delta t=0.1s$ and $N=100$ as
\begin{subequations}
  \begin{align}
    \label{equ:pend:qref}
    q_{ref} = \Bigg\{&q_{ref}(k) =
    \begin{cases}
      0 & \text{if } k \in [0,N/2) \\
      \pi & \text{if } k \in [N/2, N]
    \end{cases}
    \Bigg\} \\
    p_{ref} = \{&p_{ref}(k)=0 \;|\; k=0,\ldots,N\} \\
    u_{ref} = \{&u_{ref}(k)=0 \;|\; k=0,\ldots,N\} .
  \end{align}
\end{subequations}
Thus the reference trajectory is a step-function from the pendulum's stable
equilibrium to its unstable equilibrium. It is clearly infeasible as the desired
momentum term is zero for all time while the configuration is not, and the
desired input is zero while the system moves away from an equilibrium.

\begin{figure}
  \centering
  \includegraphics[width=0.78\columnwidth]{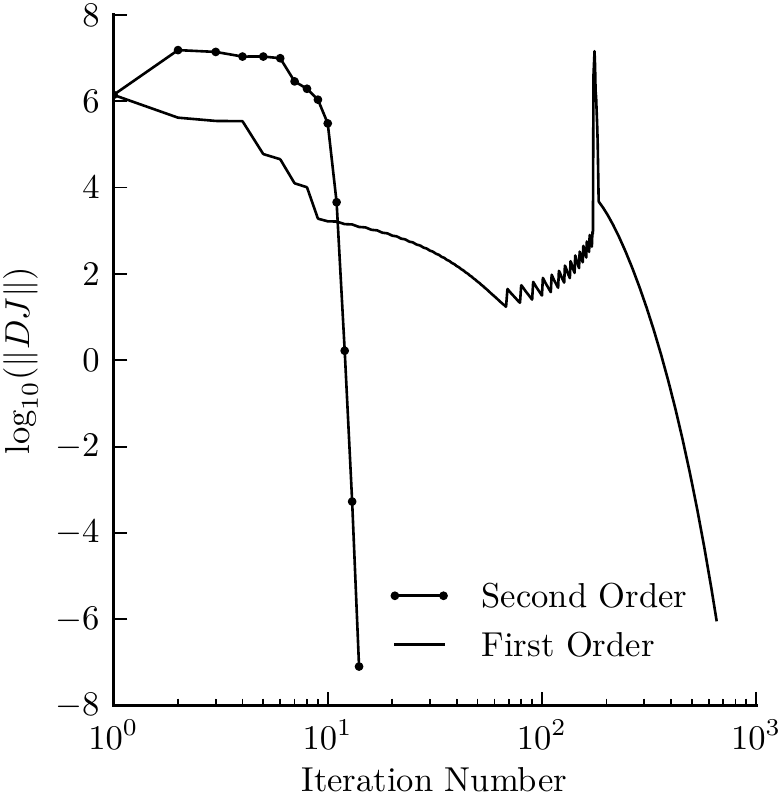}
  \caption{Convergence rate of a first-order and a second-order optimization
    method for the pendulum example.  The second order method converges to
    $10^{-6}$ in 14 iterations while the first order method converges to the
    same tolerance after 653 iterations. On an Intel i7-3770K CPU at 3.50GHz
    descent direction computation requires, on average, 28.5~ms for the
    first-order method and 55.8~ms for the second-order method resulting in the
    second-order optimization converging in approximately 4\% of the time
    required for the first-order optimization.\vspace{-12pt}}
  \label{fig:pend:convergence}
\end{figure}

An optimization routine is utilized to generate a dynamically feasible reference
trajectory where the cost function includes a weighted running cost on state and
input error as well as a weighted error of the final state.  This optimization
is performed using both a first-order method and a second-order method.  The
convergence of these two optimizations can be seen in
Fig. \ref{fig:pend:convergence} which illustrates the vastly increased
convergence rate of the second-order method.  Both methods converge to the same
feasible trajectory. It is important to note that the second-order method
requires the second derivatives presented in Sec. \ref{sec:second-deriv}.

\begin{figure}[ht!]
  \centering
  \includegraphics[width=0.78\columnwidth]{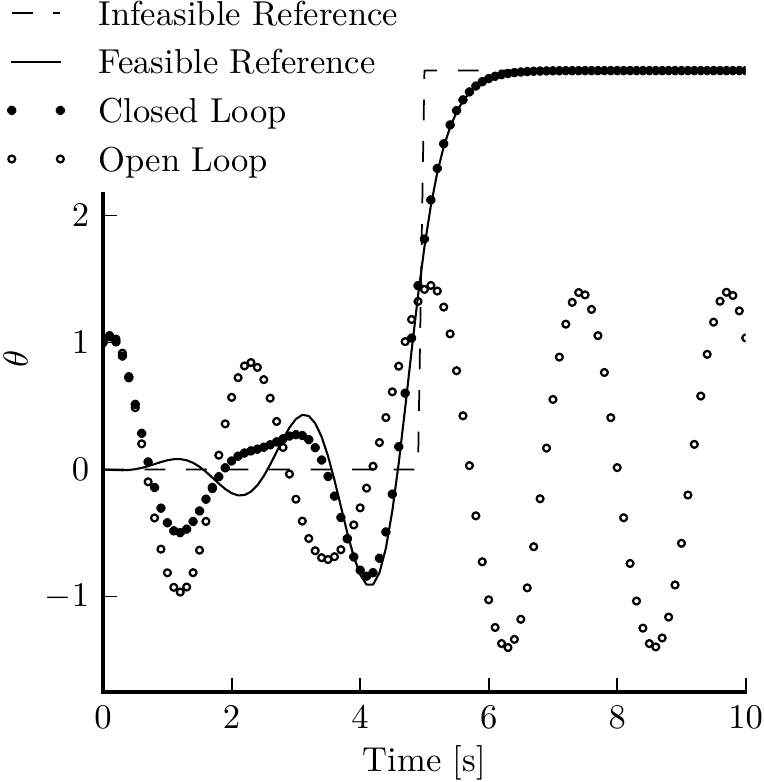}
  \caption{Figure showing the performance of the LQR regulator for the pendulum
    example.  The Infeasible Reference is the $q_{ref}$ of
    \eqref{equ:pend:qref}, the Feasible Reference is the desired trajectory
    found by the second-order optimization routine, the Closed Loop trajectory
    is produced by simulating the perturbed initial condition using
    \eqref{equ:cl:pend}, and the Open Loop trajectory uses the perturbed initial
    conditions with no stabilizing feedback.\vspace{-12pt}}
  \label{fig:pend:clsim}
\end{figure}

Once the optimization is complete, we have a dynamically feasible desired
trajectory given by $x_{d} = \{x_{d}(k)\}_{k=0}^{N}$ and $u_{d} =
\{u_{d}(k)\}_{k=0}^{N-1}$.  We now linearize the system about the
desired trajectory $(x_{d}, u_{d})$ using the derivatives of
Sec. \ref{sec:first-deriv} and then solve a Linear Quadratic Regulator (LQR)
problem to find a controller that stabilizes the system about this desired
trajectory \cite{AndersonMoore}. The LQR problem for a discrete nonlinear system
that has been linearized about a desired trajectory $(x(k),u(k))$ seeks to find
a control input $\mu(k)$ that minimizes the quadratic cost
\begin{align*}
  V(z_{k_{0}}, \mu(\cdot), k_{0}) =& \sum\limits_{k=k_{0}}^{k_{f}-1}\left[
    z^{T} (k)Q(k)z(k) + \mu^{T} (k)R(k)\mu(k) \right] \\ 
  & + z^{T}(k_{f})Q(k_{f})z(k_{f})
\end{align*}
where
\begin{align*}
  R(k) &= R^{T}(k) \geq 0 \;\;\forall \;k \in \{k_{0} \ldots (k_{f}-1)\} \\
  Q(k) &= Q^{T}(k) \geq 0 \;\;\forall \;k \in \{k_{0} \ldots k_{f}\} \\
  z_{k_{0}} &= z_{0} \\
  z(k+1) &= A(k)z(k) + B(k)\mu(k)
\end{align*}
where $z(k)$ and $\mu(k)$ are perturbations from the desired trajectory
\cite{AndersonMoore}.  The linearizations are $A(k) = \tderiv[f(x(k),
u(k))]{x(k)}$ and $B(k) = \tderiv[f(x(k), u(k))]{u(k)}$. The solution to this
discrete LQR problem is found by solving the discrete Ricatti equation:
\begin{subequations}
  \begin{align}
      \label{equ:discrete-ricatti-ofk-1}
      &P(k) = Q(k) + A^T(k) P(k+1) A(k) - \\
      &\hspace{.2in}A^T(k) P(k+1) B(k) \big[R(k) + B^T(k) P(k+1) B(k) \big]^{-1}  \notag\\
      &\hspace{.4in}B^T(k) P(k+1) A(k) \notag\\
      \label{equ:discrete-ricatti-ofk-2}
      &P_{k_f} = Q_{k_f} .
  \end{align}
\end{subequations}

The Ricatti equation is solved to find $P(k)$ by recursively evaluating
\eqref{equ:discrete-ricatti-ofk-1} backwards in time from the boundary condition
\eqref{equ:discrete-ricatti-ofk-2}.  The solution is used to calculate a 
feedback law:
\begin{equation*}
  \K(k) = \left[R(k) + B^T(k) P(k+1) B(k) \right]^{-1} B^T(k) P(k+1) B(k).
\end{equation*}
Using this feedback law with the original desired trajectory $(x_{d},u_{d})$
yields the closed-loop system 
\begin{align}
  \label{equ:cl:pend}
  &x(0) = x_{0} \nonumber \\
  &x(k+1) = f(x(k), \bar{u}(k)) \\
  &\bar{u}(k) = u_{d}(k) - \K(k)\left(x(k) - x_{d}(k)\right). \nonumber
\end{align}

To illustrate the stabilization of the controller obtained by solving the LQR
problem we perturb the initial condition of the original optimization $x_{d}(0)
= [0\;\;0]^{T}$ to give an initial condition of $x(0) = [1\;\;1]^{T}$ and then
simulate the closed loop system of \eqref{equ:cl:pend} with the perturbed
initial condition.  We also simulate the perturbed initial condition with no
feedback. The results are shown in Fig. \ref{fig:pend:clsim}, and it can be seen
that the feedback quickly stabilizes the closed-loop system to the feasible
reference.

The optimization for generating a feasible trajectory using {\tt trep} can be
found at
\texttt{\url{https://trep.googlecode.com/git/examples/papers/tase2012/pend-optimization.py}}
and the simulations of the perturbed initial conditions can be found at
\texttt{\url{https://trep.googlecode.com/git/examples/papers/tase2012/pend-closed-loop.py}}.


\subsection{Example: Marionette \label{sec:marionette}}
\begin{figure}[h!]
  \centering
  \includegraphics[height=2.75in]{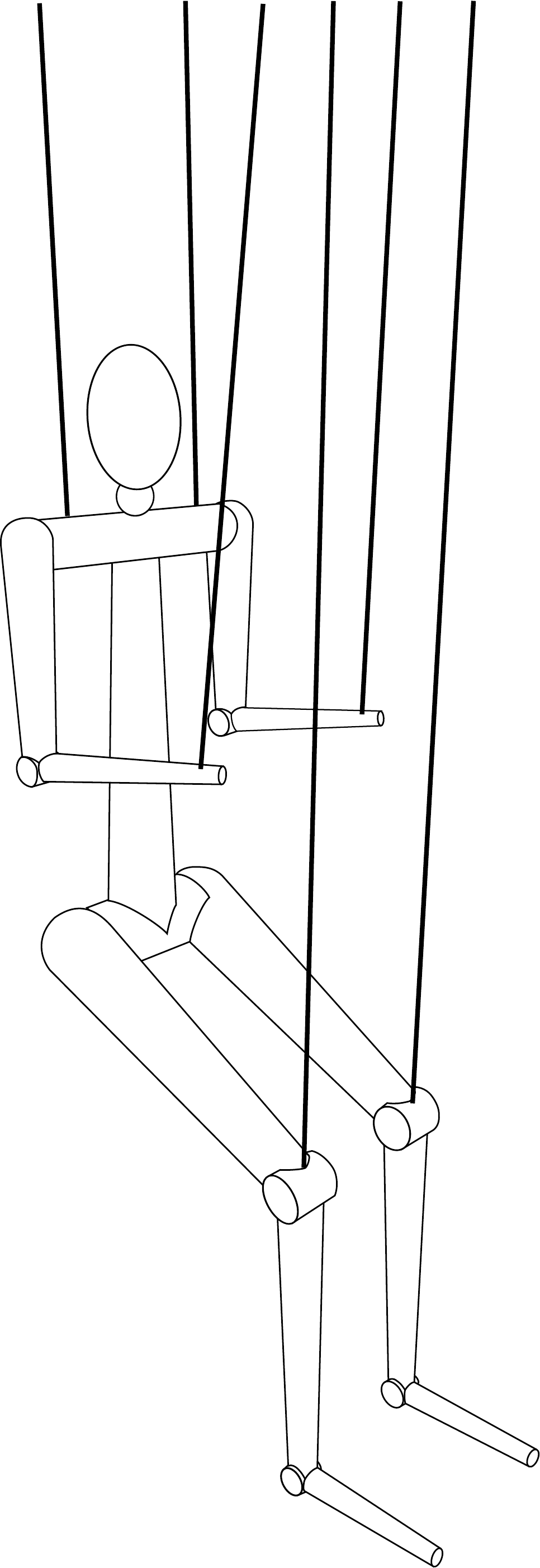}
  \caption{The marionette model has 40 configuration variables and 6
    holonomic constraints.  }
  \label{fig:marionette}
\end{figure}

As more complex example, we again use the Linear Quadratic Regulator (LQR)
method to generate a stabilizing feedback controller for the mechanical
marionette in Fig.~\ref{fig:marionette}.  The marionette has 22 dynamic
configuration variables, 18 kinematic configuration variables
\cite{ERJohnson2007}, and 6 holonomic constraints.  The corresponding
state-space model has 80 state and 18 input variables.

The marionette was simulated and linearized about a 10.0 second trajectory using
the midpoint variational integrator in {\tt trep}.  The reference trajectory was
generated by changing the string lengths of the arms and legs using $\pm 0.1
\sin(0.6 \pi t)$ input signals.  The linearization was used to create a locally
stabilizing controller by solving the discrete LQR problem with diagonal
matrices for each cost matrix with an entry of 100 for the configuration
variables and identity everywhere else. A perturbation of $0.1$ $rad$ was then
added to the initial condition of the vertical orientation of the torso and the
simulation was performed with and without the added stabilizing feedback
controller.

\begin{figure}
  \centering
  \includegraphics[height=2.65in]{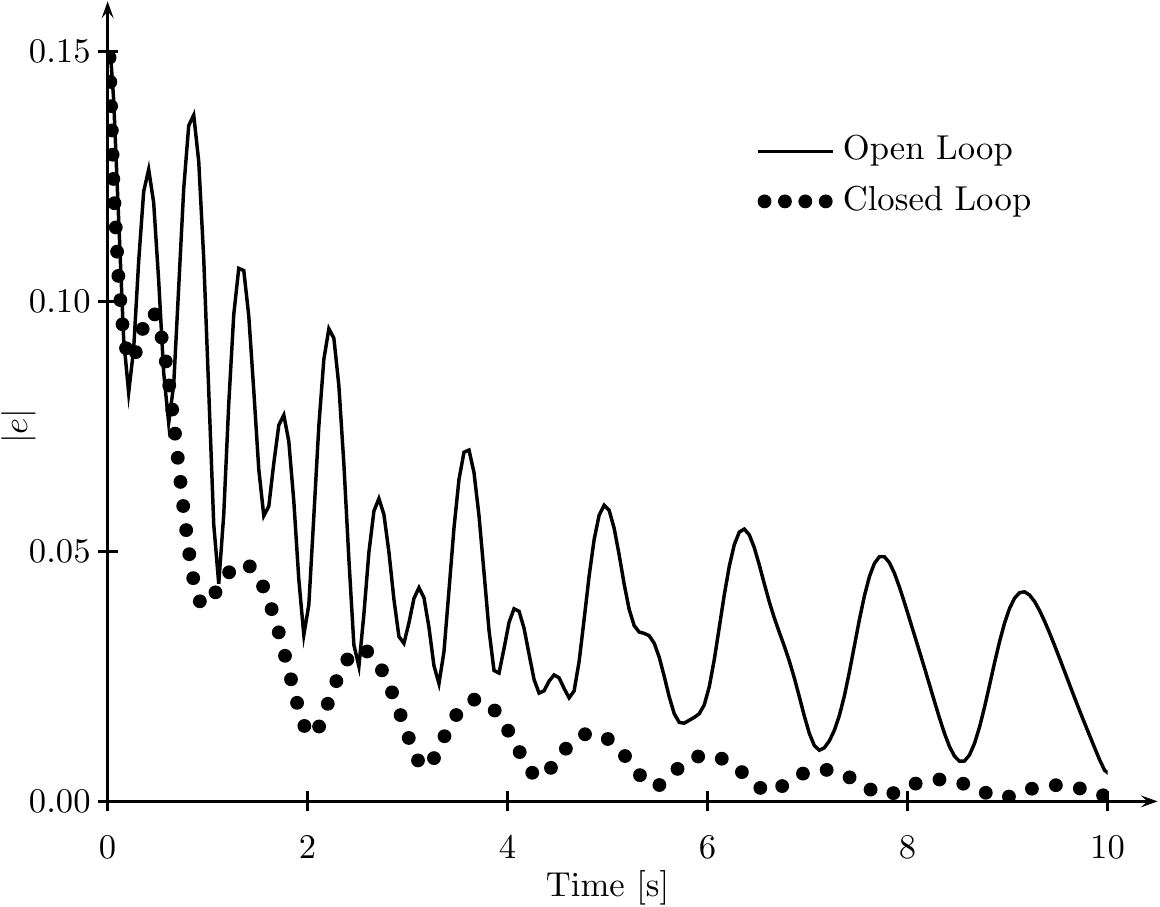}
  \caption{The discrete LQR feedback law significantly improves the
    norm of the error response of the marionette compared to the open-loop
    simulation. 
  }
  \label{fig:ex1-error}
\end{figure}

\begin{figure}
  \centering
  \includegraphics[height=2.65in]{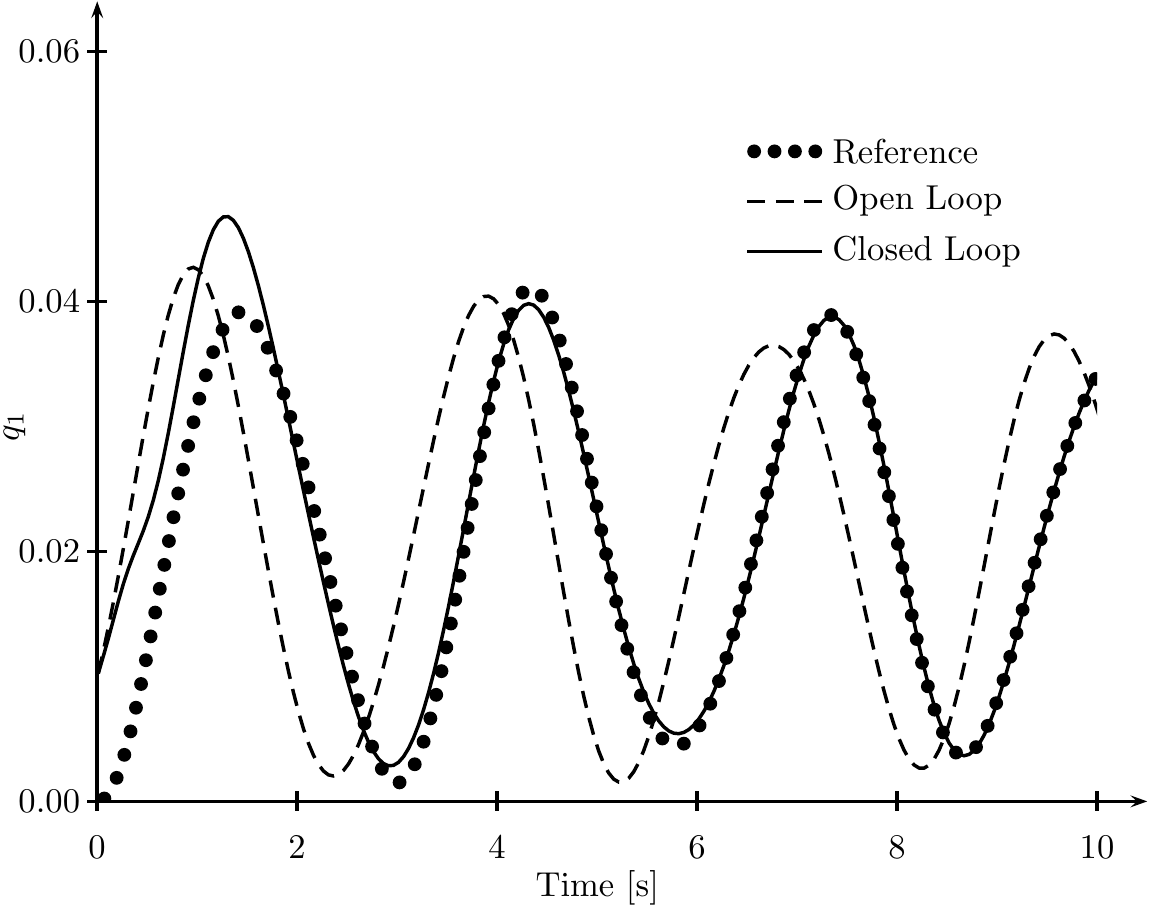}
  \caption{The discrete LQR feedback law also significantly improves the
    individual 
    error response of the marionette compared to the open-loop
    simulation.  This is the trajectory of the  configuration of the
    vertical orientation of the torso.  }
  \label{fig:ex1-trajerror}
\end{figure}

The norm of the resulting error between the perturbed and original trajectories
is shown in Fig.~\ref{fig:ex1-error} and the trajectory of the vertical
orientation of the torso is shown in Fig.~\ref{fig:ex1-trajerror} as an example
of stabilization of one of the states.  As expected, the closed-loop trajectory
quickly converges to the reference trajectory while the errors in the open-loop
trajectory persist throughout the time horizon.  The ability to generate locally
stabilizing feedback laws for complex systems that are simulated with
variational integrators is a useful application of the methods described
here. The source code for this example can be obtained at
\texttt{\url{https://trep.googlecode.com/git/examples/papers/tase2012/marionette.py}}.

The optimization was performed on an Intel i7-2760QM CPU at 2.40GHz.  The
simulation takes approximately 2.11 ms per step, the linearization takes
approximately 1.12 ms per step.  The second-order linearization takes
approximately 22.15 ms per step, though it was not required for this example.


\section{Conclusion\label{sec:conclusion}}

Variational integrators are an appealing alternative to numerically integrating
continuous equations of motion.  By representing variational integrators as
discrete dynamic systems and calculating the linearization of the associated
one-step map, their utility is extended to applications requiring analysis and
optimal control.  This approach reduces complexity, potential for error, and
extraneous work compared to using a variational integrator for simulation while
doing the analysis and optimization in the continuous domain with a separate set
of equations.  Moreover, it leads to feedback laws that are expressed purely in
terms of configuration variables (instead of configurations and configuration
velocities).

The methods described here can be efficiently implemented by using a recursive
tree representation to calculate the required derivatives of the discrete
Lagrangian and forcing function.  The approach accommodates external forcing and
holonomic constraints as described here, and is compatible with kinematic
configuration variables \cite{ERJohnson2007,JohnsonThesis}.  Additionally, this
method could be extended to calculate higher derivatives if needed.

The authors have additionally used this framework to implement projection
operator-based trajectory optimization \cite{HauserSaccon} using variational
integrators as the representation of the dynamics (these results will be
presented in a future article).

\section*{Acknowledgment}
This material is based upon work supported by the National Science Foundation
under award CNS 1329891 and IIS-0757378.  Any opinions, findings, and
conclusions or recommendations expressed in this material are those of the
authors and do not necessarily reflect the views of the National Science
Foundation.

\ifCLASSOPTIONcaptionsoff
  \newpage
\fi



\bibliographystyle{IEEEtran}
\bibliography{bibliography}
\end{document}